\documentclass{amsart}
\usepackage{mathrsfs}
\usepackage{amssymb, amscd}
\usepackage[matrix, arrow]{xy}

\newtheorem{theorem}{Theorem}[section]
\newtheorem{lemma}[theorem]{Lemma}
\newtheorem{corollary}[theorem]{Corollary}
\newtheorem{proposition}[theorem]{Proposition}
\newtheorem{definition}[theorem]{Definition}

\DeclareMathOperator{\Spec}{Sp}
\DeclareMathOperator{\Hom}{Hom}
\DeclareMathOperator{\RHom}{RHom}

\newcommand{\RHomC}{\RHom_{\C}}

\DeclareMathOperator{\End}{End}

\DeclareMathOperator{\Ho}{Ho}
\DeclareMathOperator{\id}{id}
\DeclareMathOperator{\op}{op}

\DeclareMathOperator{\ev}{Ev}

\newcommand{\mQ}{{\mathbb Q}}
\newcommand{\mZ}{{\mathbb Q}}

\newcommand{\bT}{{\mathbb T}}
\newcommand{\bI}{{\mathbb I}}
\newcommand{\bQ}{{\mathbb Q}}
\newcommand{\bZ}{{\mathbb Z}}

\newcommand{\C}{{\mathcal C}}
\newcommand{\D}{{\mathcal D}}
\newcommand{\E}{{\mathcal E}}

\newcommand{\iso}{\cong}

\newcommand{\sm}{\wedge}
\newcommand{\ssm}{\ltimes}

\newcommand{\boxprod}{\mathbin{\square }}

\newcommand{\wh}{\widehat}
\newcommand{\mc}{\colon \,}

\renewcommand{\to}{\longrightarrow}
\newcommand{\varrow}[1]{\hbox to #1{\rightarrowfill}}
\newcommand{\varl}[2]{\stackrel{#2}{\hbox to #1{\leftarrowfill}}}
\newcommand{\varrx}[1]{\stackrel{#1}{\hbox to {1cm}{\rightarrowfill}}}
\newcommand{\varr}[1]{\stackrel{#1}{\longrightarrow}}

\newcommand{\ch}{Ch_{\bQ}}
\newcommand{\EG}{\cE(\cG)}
\newcommand{\cB}{\mathcal{B}}
\newcommand{\B}{\mathcal{B}}
\newcommand{\Ba}{\mathcal{B}_a}
\newcommand{\Bt}{\mathcal{B}_t}

\newcommand{\cG}{\mathcal{G}}
\newcommand{\cF}{\mathcal{F}}
\newcommand{\cO}{\mathcal{O}}
\newcommand{\cA}{\mathcal{A}}
\newcommand{\cE}{\mathcal{E}}
\newcommand{\cD}{\mathcal{D}}
\newcommand{\cT}{\mathcal{T}}
\newcommand{\cR}{\mathcal{R}}
\newcommand{\cS}{\mathcal{S}}
\newcommand{\tg}{\widetilde{g}}
\newcommand{\bb}{s}
\newcommand{\GrO}{\cO_{\cF *}}
\newcommand{\BC}{\mathcal{B}_a}
\newcommand{\BCt}{\mathcal{B}_t}
\newcommand{\EBCa}{\cE(\cB_a)}
\newcommand{\EBCo}{\cE(\cB_a)\langle 0 \rangle}
\newcommand{\EBCto}{\cE(\cB_t)\langle 0 \rangle}
\newcommand{\EBCt}{\cE(\cB_t)}
\newcommand{\HomC}{\Hom_{\C}}
\newcommand{\HomR}{\Hom_{\cR}}
\newcommand{\CHom}{\Hom_{\C}}

\newcommand{\Mod}{\mathcal{M}od\mbox{-}}
\newcommand{\Et}{\cE_t}
\newcommand{\Ea}{\cE_a}

\newcommand{\Qeq}{\simeq_Q}
\newcommand{\Prod}{\Pi}
\newcommand{\Osum}{\bigoplus}
\newcommand{\osum}{\oplus}
\newcommand{\dga}{dg\cA}

\begin{document}

\title{An algebraic model for rational $S^1$-equivariant stable homotopy
theory}
\date{\today; 2000 AMS Math.\ Subj.\ Class.: 55P62, 55P91, 55P42, 55N91, 18E30}
\author{Brooke Shipley}
\thanks{Research partially supported by an NSF grant.}
\address{Department of Mathematics \\ Purdue University \\
West Lafayette, IN 47907 
\\ USA}
\email{bshipley@math.purdue.edu}

\begin{abstract}
Greenlees defined an abelian category $\cA$ whose derived category is 
equivalent to the rational $S^1$-equivariant stable homotopy category
whose objects represent rational $S^1$-equivariant cohomology theories.
We show that in fact the model category of differential graded objects
in $\cA$ models the whole rational $S^1$-equivariant stable homotopy theory.
That is, we show that there is a Quillen equivalence between $\dga$ and
the model category of rational $S^1$-equivariant spectra, before the
quasi-isomorphisms or stable equivalences have been inverted.
This implies that all of the higher order structures such as mapping spaces,
function spectra and homotopy (co)limits are reflected in the 
algebraic model.  The construction of this equivalence involves calculations
with Massey products.  In an appendix we show that Toda brackets, and
hence also Massey products, are determined by the derived category.    
\end{abstract}

\maketitle

\section{Introduction} 

In~\cite{Greenlees}, Greenlees defined an abelian category $\cA$ and
showed that its derived category is equivalent to the rational
$\bT$-equivariant stable homotopy category where $\bT$ is the circle
group.  We strengthen this result by showing that this derived
equivalence can be lifted to an equivalence 
on the underlying model categories, before the quasi-isomorphisms or
stable equivalences have been inverted.

\begin{theorem}\label{thm-main-intro}
The model category of rational $\bT$-equivariant spectra is Quillen
equivalent to the model category of differential graded objects in $\cA$.
\end{theorem}
\noindent The definition of $\cA$ and the model category
on differential graded objects in $\cA$, $\dga$, are recalled in
Section~\ref{sec-model}.

As mentioned in~\cite[16.1]{Greenlees}, one of the applications of this
stronger equivalence is that the algebraic models for the smash product
and function spectra given in~\cite[Part IV]{Greenlees} are natural.
This also shows that the higher order structures such as mapping spaces
and homotopy (co)limits are captured by the algebraic 
model~\cite{DK,dwyer-spalinski}.  Since
$\cA$ has injective dimension one, calculations of these higher order
structures are quite practical, as demonstrated here and 
in~\cite{Greenlees}.  

\cite[Part III]{Greenlees} discusses many motivations for
studying rational $\bT$-equivariant spectra, including its connection
to algebraic $K$-theory, equivariant topological $K$-theory and Tate
cohomology.  This algebraic model has also been used to define a
model for rational $\bT$-equivariant elliptic cohomology~\cite{GHR}.
Another reason for studying rational $\bT$-equivariant spectra
is that the circle is the simplest infinite compact Lie group.
The arguments used in~\cite[5.6.1]{Greenlees} to show that
the derived category, $\cD(\cA)$, is equivalent to the homotopy
category of rational $\bT$-equivariant spectra, $\Ho(\bT\mbox{-spectra})$,
rely heavily on the fact that the circle is rank one.  The general
approach here though, using Quillen model categories and
Morita equivalences, does apply to compact Lie groups of higher 
rank.  Algebraic models of rational $T^r$-equivariant spectra, for
tori of any rank, are considered in~\cite{GS}.  

The main new ingredient here is a Morita equivalence for stable model
categories considered in~\cite{ss3}, see Section~\ref{sec-Morita}.
In general, this Morita equivalence models a stable model category
with a set of small generators by modules over a ring spectrum
(with many objects).  For rational stable model categories, this
can be simplified further as (unbounded) modules over a rational
{\em differential graded ring with many objects}; see 
Definition~\ref{def-ringoids} or~\cite[1.3]{ssme}.  
The following statement is proved as Corollary~\ref{cor-Gabriel-t} below.
Definition~\ref{def-Quillen} recalls the notion of a Quillen 
equivalence.

\begin{theorem}\label{thm-Gabriel-intro}
The model category of rational $\bT$-equivariant spectra is
Quillen equivalent to the model category of (right) modules over
a rational differential graded ring $\E(\Bt)$ with one object
for each closed subgroup in $\bT$.
\[
\bT\mbox{-spectra} \Qeq \Mod\E(\Bt)
\]
\end{theorem}

This gives an algebraic model for rational $\bT$-equivariant spectra,
but $\E(\Bt)$ is large and difficult to make explicit.  Hence
this algebraic model is not very practical.  Instead, it is useful
as a stepping stone to create a Quillen equivalence between
$\bT$-spectra and differential graded objects in $\cA$.  As another
stepping stone we apply this same process to $\dga$.  
The model category of differential graded objects in $\cA$
is Quillen equivalent to (right) modules over a rational differential
graded ring $\E(\Ba)$ with one object for each closed subgroup in $\bT$;
see Corollary~\ref{cor-Gabriel} below.  We then show that these
two differential graded rings (with many objects) $\E(\Bt)$ and $\E(\Ba)$
are quasi-isomorphic.   
Since quasi-isomorphisms of differential
graded rings (with many objects) induce Quillen equivalences of the 
associated categories of modules by~\cite[4.3]{ss1}, see 
also~\cite[A.1.1]{ss3} or~\cite[A.1]{ssme},
Theorem~\ref{thm-main-intro} follows given a zig-zag of
quasi-isomorphisms between $\E(\Ba)$ and $\E(\Bt)$.  This is carried
out in Sections~\ref{sec-endo},~\ref{sec-end-t} and~\ref{sec-map}.

\begin{theorem}\label{thm-zig-intro}
The model category of differential graded objects in $\cA$ is 
Quillen equivalent to modules over the differential graded ring $\E(\Ba)$. 
Also, there is a zig-zag of quasi-isomorphisms 
between $\E(\Ba)$ and $\E(\Bt)$
which induces a zig-zag of Quillen equivalences between the
associated model categories.  Hence we have a chain of
Quillen equivalences
\[ 
\bT\mbox{-spectra}\Qeq \Mod\E(\Bt) \Qeq \Mod\E(\Ba) \Qeq \dga.
\]
\end{theorem}

In more detail, in Section~\ref{sec-endo} we explicitly describe
$\E(\Ba)$ and a quasi-isomorphic sub-ring $\Ea$.  In 
Section~\ref{sec-end-t}, we modify $\E(\Bt)$ to define a more amenable
ring $\Et$.  Then in Section~\ref{sec-map} we construct an
intermediary differential graded ring $\cS$ and two quasi-isomorphisms
$\Ea \xleftarrow{\eta} \cS \xrightarrow{\varphi} \Et$.  
This construction involves Massey products
for differential graded rings with many objects which are defined in 
Definition~\ref{def-Massey}.

In particular, the proof of Theorem~\ref{thm-zig-intro} uses the fact that
the higher order products for $\Ea$ and $\Et$ agree. These products
agree because there are Quillen equivalences
\[ \Mod\Ea \Qeq \Mod\E(\Ba) \Qeq \dga \]
and
\[ \Mod\Et \Qeq \Mod\E(\Bt) \Qeq \bT\mbox{-spectra} \]
and $\Ho(\bT\mbox{-spectra})$ and $\cD(\cA)$ are equivalent as
triangulated categories by~\cite[5.6.1]{Greenlees}.  This invariance of
higher order products does not seem to appear in the literature, so we
include it in an appendix on Toda brackets. 
A related statement, that quasi-isomorphisms preserve
Massey products, appears in~\cite[1.5]{may}.
Specifically, we prove the following as Theorem~\ref{thm-tb-eq}.

\begin{theorem} If $\varphi\mc \cT \to \cT'$ is an exact equivalence of
triangulated categories then the Toda brackets for $\cT$ and $\cT'$
agree.
\end{theorem}

{\em Acknowledgments:} I thank Bill Dwyer for a discussion
which sparked my interest in this project. I also thank John 
Greenlees and
Jeff Smith for many conversations throughout this project.  Finally, I
thank Mike Mandell and Peter May for sharing their
equivariant expertise and Mark Mahowald for sharing his understanding 
of Toda brackets.

\section{Algebraic Model}\label{sec-model}

In this section we recall Greenlees' algebraic model for rational
$\bT$-equivariant  spectra; for more detail 
see~\cite[Chapters 4, 5; Appendix E.1]{Greenlees}.  We then show
that this model is a $\ch$-model category; see Definition~\ref{def-chain model
cat}.  Finally we recall certain objects called the
algebraic basic cells.

First we recall some preliminary definitions.

\begin{definition}~\cite[4.5, 4.6, 5.2]{Greenlees}
{\em Let $\cF$ be the set of finite subgroups of $\bT$.  Let $\cO_{\cF}$  
be the ring of operations, $\Prod_{H\in \cF} \bQ[c_H]$ with $c_H$ in degree 
$-2$.  Let $e_H$ be the idempotent which is projection onto the
factor corresponding to $H$.  Let $c$ be the {\em total Chern class} with
$e_H c = c_H$.   Let $\cE$ be the multiplicative set of all {\em Euler classes},
$\{c^v\  |\  v \mc \cF \to \bZ_{\geq 0} \mbox{ with finite support}\}$, where
$e_Hc^v = c_H^{v(H)}$.  Let \ $t_*^{\cF}= \cE^{-1}\cO_{\cF}$.  As a vector
space $(t_*^{\cF})_{2n}$ is $\Prod_{H\in \cF} \bQ$ for $n \leq 0$ and is 
$\Osum_{H \in\cF}\bQ$
for $n > 0$.  
}\end{definition}

We now define the abelian category $\cA$.  The category of differential
graded objects in $\cA$ is then the standard algebraic model of rational 
$\bT$-equivariant spectra.

\begin{definition}~\cite[5.4.1]{Greenlees}
{\em The objects of $\cA$ are maps $N \varr{\beta} t_*^{\cF} \otimes V$ of
$\cO_{\cF}$-modules with $V$ a graded vector space such that $\cE^{-1}\beta$
is an isomorphism.  $N$ is referred to as the {\em nub} and $V$ is the 
{\em vertex}. 
Morphisms are commutative squares
\[\begin{CD}
M @>{\theta}>> N \\
@V{\beta}VV  @VV{\beta'}V \\
t_*^{\cF} \otimes U @>>{1 \otimes \phi}> t_*^{\cF} \otimes V.
\end{CD}\]
}\end{definition}

The condition that $\cE^{-1}\beta$ is an isomorphism is equivalent
to requiring that the kernel $K$ and cokernel $C$ of $\beta$ are
$\cF$-finite torsion modules~\cite[4.5.1, 4.5.2, 4.6.6]{Greenlees}.  That is, 
they are torsion modules with respect to the total Chern class and they 
decompose
as direct sums, $K = \Osum_H e_H K$ and $C = \Osum_H e_H C$. 

Greenlees, in~\cite[5.6.1]{Greenlees}, shows that the derived category
of differential graded objects in $\cA$ is equivalent to the homotopy
category of rational $\bT$-equivariant spectra.  Here we recall the
model category structure on $\dga$.  The associated homotopy category
is equivalent to $\cD(\cA)$.
 
\begin{proposition}~\cite[Appendix B]{Greenlees}
The category $dg\cA$ of differential graded objects in $\cA$ with
cofibrations the monomorphisms, weak equivalences the quasi-isomorphisms,
and fibrations determined by the right lifting property is
a model category. 
\end{proposition}

The model category $dg\cA$ is also tensored, cotensored and enriched
over rational chain complexes in a way that is compatible with
the model category structure.   We call such a model category 
a $\ch$-model category.  This
is the analogue of a simplicial model category~\cite[II.2]{Q} with simplicial
sets replaced by rational chain complexes.  See~\cite[Ch. 4]{hovey} for
a more general definition. Here $\ch$ is the category 
of unbounded rational chain complexes 
with the projective model category structure~\cite[2.3.3]{hovey} 
and the standard closed symmetric monoidal structure~\cite[4.2.13]
{hovey}.  For definitions of {\em tensor, cotensor} and {\em enriched}
see~\cite[1.2, 3.7]{kelly}.

\begin{definition} \label{def-chain model cat}
{\em A {\em $\ch$-model category} is a complete and cocomplete model category 
$\C$ which is
tensored and enriched (denoted $\HomC$) over the category
$\ch$, has cotensors with finite dimensional complexes in $\ch$ and 
satisfies the following compatibility axiom (CM):\\
(CM) For every cofibration $A\varr{i} B$ and every fibration $X\varr{p} Y$ in 
$\C$ the induced map
\[ \HomC(i, p) \mc \HomC(B,X) \ \varrow{1cm} \ \HomC(A,X) \times_{\HomC(A,Y)}
\HomC(B,Y) \]
is a fibration in $\ch$. If in addition one of the maps
$i$ or $p$ is a weak equivalence, then $\HomC(i,p)$ is a trivial fibration.
We use the notation $K\otimes X$ and $X^K$ to denote the tensors and cotensors
for $X$ in $\C$ and $K$ a chain complex.}
\end{definition}

For example, $\ch$ is itself a $\ch$-model category~\cite[4.2.13]{hovey}
as is the projective model category defined 
in~\cite[A.1.1]{ss3} or~\cite[A.1]{ssme} of unbounded differential graded
modules over any rational differential graded ring (with many objects);
see Definition~\ref{def-ringoids}.  

\begin{proposition}\label{prop-ch-mc}
The model category on $dg\cA$ is a $\ch$-model category.
\end{proposition}

\begin{proof}
For $K \in \ch$ and $A \in dg\cA$ define $K \otimes A$ in $dg\cA$ by
\[K \otimes A = \Osum_n \Sigma^n (K_n \otimes_{\bQ} A)\]
with differential $d(k_n \otimes a) = dk_n \otimes a
+ (-1)^n k_n \otimes da$.  Let $\bQ[n]$ denote the rational chain complex
with $\bQ$ in degree $n$.  Then for $A, B$ in $dg\cA$ define 
$\CHom(A,B)$ in $\ch$ by $\CHom(A, B)_n = \cA(\bQ[n] \otimes A, B)$ where
$\cA(-,-)$ is the rational vector space of maps of underlying graded
objects in $\cA$ (ignoring the differential).  The differential
for $\CHom(A,B)$ is given by $df_n = d_B f_n + (-1)^{n+1} f_n d_A$.    
For $K$ a finite dimensional complex, define $A^K$ as $\Hom_{\ch}(K,\mQ[0]) 
\otimes A$. This cotensor is the right adjoint of the tensor $K \otimes -$.   
So we have the following isomorphisms (where defined)
\[ \ch(K, \HomC(A,B)) \iso dg\cA( K \otimes A, B) \iso dg\cA(A, B^K). \]   

As in~\cite[II.2 SM7(b)]{Q} or~\cite[4.2.1, 4.2.2]{hovey},  the compatibility 
axiom in Definition~\ref{def-chain model cat} has an equivalent adjoint form involving the tensor.
For $f\mc K \to L$ in $\ch$  and $g\mc A \to B$ in $dg\cA$, define
$f \boxprod g \mc K \otimes B \amalg_{K \otimes A} L \otimes A \to L \otimes B$.
Then Axiom (CM) from Definition~\ref{def-chain model cat} is equivalent to
requiring that
if $f$ and $g$ are cofibrations, then 
$f \boxprod g$ is a cofibration which is trivial if $f$ or $g$ is.      
Since $\ch$ is a cofibrantly generated model category we only need to 
check this
property when $f$ is one of the generating cofibrations or generating
trivial cofibrations~\cite[3.5]{ss1}.

Let $D^n$ be the acyclic rational chain complex with $\bQ$ in degrees $n$ 
and $n-1$.  Then the generating cofibrations of $\ch$ are the inclusions 
$i_n \mc \bQ[n-1] \to D^n$ for $n$ an integer and the generating trivial 
cofibrations are the maps $j_n \mc 0 \to D^n$ for $n$ an 
integer~\cite[2.3.3]{hovey}.      
The source of the map $i_n \otimes g$ is the cofiber of $g$ and the
target is $D^n \otimes B$.   So for $g$ a monomorphism,  
$i_n \boxprod g$ is a monomorphism because it is $g$ on one summand and 
$\id_B$ on the other.   Since $D^n \otimes B$ is always acyclic and the cofiber
of a quasi-isomorphism is acyclic, if $g$ is a trivial cofibration, then  
$i_n \boxprod g$ is a quasi-isomorphism.    Also, $j_n \boxprod g \mc 
D^n \otimes A \to D^n \otimes B$ is a quasi-isomorphism between acyclic 
complexes for any monomorphism $g$. 
\end{proof}

For the Morita equivalence discussed in the next section we need
to define the algebraic cells which can be used to build any object
in $\dga$. 
The natural building blocks for $\bT$-equivariant spectra are
the $G$-cells $\bT_+ \sm_H S^n$.   Rationally, these can be simplified by
using idempotents in the Burnside
ring for a finite group, $A(H)= [S^0, S^0]^H \iso \Prod_{K\leq H} \bQ$.
Following Greenlees, we suppress the rationalization in our notation
for spectra.  
Let $e_K$ be the idempotent which is projection onto the $K$th factor.

\begin{definition}~\cite[2.1.2]{Greenlees}
{\em The {\em geometric basic cells} are the rational $\bT$ spectra 
$\sigma^0_{\bT} = S^0$ and $\sigma^0_H = \bT_+ \sm_H e_H S^0$ 
for $H$ a finite subgroup of $\bT$.
}\end{definition}

Via the equivalences between $\cD(\cA)$ and rational $\bT$-equivariant
spectra, these geometric basic cells correspond to the following algebraic
basic cells.

\begin{definition}~\cite[5.8.1]{Greenlees}
{\em For $H$ a finite subgroup, let $\bQ(H)$ denote the
$\bQ[c_H]$-module $\bQ$ regarded as an $\cO_{\cF}$-module.  Then
the {\em algebraic basic cells} are the objects $L_H =(\Sigma \bQ(H) \to 0)$  and $L_{\bT} = (\cO_{\cF} \to t_*^{\cF})$ in $\cA$.} 
\end{definition}

We need fibrant replacements of these algebraic 
basic cells which are inclusions 
into their respective injective hulls.  The injective dimension of $\cA$
is one, so the injective hulls are simple to construct.  First we define 
the basic injective objects.   

\begin{definition}~\cite[2.4.3, 3.4.1, 5.2.1, 5.5.1]{Greenlees}
{\em Let $\bI(H) = \Sigma^{-2} \bQ[c_H, c_H^{-1}]/\bQ[c_H]$.  Let
$\bI = \Osum_H \bI_H  \iso \Sigma^{-2} t_*^{\cF} /\cO_{\cF}$.  The 
underlying vector spaces are $\bI(H)_{2n} = \bQ$ for $n \geq 0$ and
$\bI_{2n} = \Osum_H \bQ$ for $n \geq 0$.  
Then $(\bI \to 0)$, $ (\bI(H) \to 0)$,  and 
$(t_*^{\cF} \otimes V \to 
t_*^{\cF}\otimes V)$ for any graded vector space $V$ are injective objects
in $\cA$.  
}\end{definition}
   
The construction of fibrant approximations in~\cite[Appendix B]{Greenlees} 
produces injective hulls for $L_H$ and $L_{\bT}$ which we denote by
$I_H$ and $I_{\bT}$ respectively.  Let $P(I)= I \ssm \Sigma^{-1}I$ be the 
acyclic path object
in $dg\cA$ with underlying graded object $I \osum \Sigma^{-1} I$.   
For a finite subgroup $H$ of $\bT$, let $\Sigma \bI(H) \ssm \Sigma^2 
\bI(H)$ be the pullback in the following square.
\[\begin{CD}
\Sigma \bI(H) \ssm \Sigma^2 \bI(H) @>>> P(\Sigma^3 \bI (H)) \\
@VVV @VVV \\
\Sigma \bI(H) @>>> \Sigma^3 \bI(H)
\end{CD}\] 
The bottom map is zero in degree one and an isomorphism in all higher
degrees.   Let $I_H$ denote the object $(\Sigma \bI(H) \ssm \Sigma^2 \bI(H) \to 0)$ 
in $dg\cA$.  The inclusion $L_H \to I_H$ is a quasi-isomorphism and $I_H$ is
fibrant since it is built from standard injectives.

For $\bT$ itself, let $t_*^{\cF} \ssm \Sigma \bI$  be the 
pullback in the following square
\[\begin{CD}
t_*^{\cF} \ssm \Sigma \bI  @>>> P(\Sigma^2 \bI)  \\
@VVV @VVV \\
t_*^{\cF} @>>> \Sigma^2 \bI
\end{CD}\] 
The bottom map is zero in degrees $n \leq 0$ and an isomorphism in
all positive degrees.  Let $I_{\bT}$ denote the object $(t_*^{\cF} \ssm 
\Sigma \bI \to t_*^{\cF})$ in $\dga$.  The inclusion $L_{\bT} \to I_{\bT}$ is a
quasi-isomorphism and $I_{\bT}$ is fibrant since it is built from
standard injectives. 

\begin{proposition}
$I_H$ and $I_{\bT}$ are fibrant replacements of $L_{H}$ and $L_{\bT}$
respectively.
\end{proposition}

\section{The Morita equivalence}\label{sec-Morita}

Gabriel \cite{Gabriel} proved a Morita theorem which shows  
that any cocomplete abelian category with 
a set of finitely generated projective generators is a category of 
$\E$-modules for some ring $\E$ (with many objects).  See 
Definition~\ref{def-ringoids} below.  The ring $\E$ in question may be
taken to be the endomorphism ring of the projective generators; 
this is naturally a ring with one object for each generator. 
An enriched variant of this Morita equivalence was considered
in~\cite{ss3}.  Here we consider another enriched variant 
which shows that a model category with a set of small generators which is 
compatibly enriched over chain complexes is Quillen equivalent to
a category of differential graded modules over a differential graded
ring (with many objects).

\begin{definition}\label{def-ringoids}
{\em A {\em ring with objects $\cG$} is an $Ab$-category with object set
$\cG$. This is also a category enriched over abelian groups,~\cite[1.2]{kelly}
or an $Ab$-module,~\cite[4.1.6]{hovey}.  
A ring is then an $Ab$-category with one object.  Similarly, $\cR$, a 
{\em rational differential graded ring with objects $\cG$}, is
a category enriched over $\ch$, or a $\ch$-category or $\ch$-module.   
A (right) $\cR$-module is a
$\ch$-functor $M \mc \cR^{\op} \to \ch$.  
}\end{definition}

Given $\cR$, a rational differential graded ring with objects $\cG$, 
and objects $G, G', G''$ in $\cG$ there is a morphism chain complex
$\HomR(G,G') \in \ch$ and coherently associative and unital composition maps
\[ \HomR(G', G'') \otimes \HomR(G, G') \to \HomR(G,G'').\]     
An $\cR$-module $M$ is then determined by chain complexes $M(G)$ for
each object $G \in \cG$ and coherently associative and unital action
maps
\[ M(G') \otimes \HomR(G,G') \to M(G). \] 
For each object $G$ in $\cG$ there is a free $\cR$-module, $F_G^{\cR}$,
defined by $F_G^{\cR}(G')=\cR(G',G)$.  
The projective model structure on the category of $\cR$-modules, $\Mod\cR$, is 
established
in~\cite[A.1.1]{ss3} or~\cite[A.1]{ssme} using~\cite[2.3]{ss1}.  
The weak equivalences are the quasi-isomorphisms and
the fibrations are the epimorphisms.  With the standard action of $\ch$, 
$\Mod\cR$ is a $\ch$-model category. 

\begin{definition} \label{def-EP}
{\em Let $\cG$ be a set of objects in a $\ch$-model category $\C$.
We denote by $\EG$ the full $\ch$-subcategory of $\C$ with objects $\cG$,
i.e., $\EG(G,G')=\HomC (G,G')$ is a differential graded ring
with many objects, one for each element of $\cG$. We let
\[ \HomC (\cG,-) \, : \,  \C \ \varrow{1cm} \ \Mod\EG\]
denote the  functor given by $\HomC (\cG,Y)(G) = \HomC (G,Y)$.}
\end{definition}

Note that, if $\cG= \{ G \}$ has a single element,
then $\EG$ is determined by the single differential graded ring,
$\End (G) = \HomC (G,G)$.

For our enriched variant of Gabriel's theorem it is crucial that
a $\ch$-model category is a {\em stable model category}, that is,
a pointed model category where the suspension functor is a self-equivalence
on the homotopy category~\cite{ss3}.   The homotopy category
of a stable model category is triangulated~\cite[7.1.1, 
7.1.6]{hovey},~\cite{verdier}.  
The suspension and loop functors are inverse equivalences which induce
the shift functor.  The cofiber and fiber sequences agree
up to a sign and induce the triangles.

\begin{proposition}\label{prop-chain-stable}
If $\C$ is a $\ch$-model category then $\C$ is a stable model category.
If $X$ is a cofibrant object in $\C$ and $Y$ is a fibrant object in $\C$,
then there is a natural isomorphism of graded abelian groups 
$H_* \Hom_{\C}(X,Y) \iso [X,Y]_*^{\Ho(\C)}$. 
\end{proposition}

\begin{proof}
First, $\C$ is pointed since $\ch$ is pointed~\cite[4.2.19]{hovey}.
For a cofibrant object $X$ in $\C$, a cylinder object~\cite[I.1]{Q} is given
by $(D^1 \osum \mQ[0]) \otimes X$.  
It is a cylinder object because 
$X \amalg X \iso (\mQ[0] \osum \mQ[0]) 
\otimes X  \varr{i} (D^1 \osum \mQ[0]) \otimes X$ is
a cofibration, by an adjoint form of the compatibility axiom (CM), 
see~\cite[II.2, SM7(b)]{Q},
and $(D^1 \osum \mQ[0]) \otimes X \to \mQ[0] \otimes X \iso X$ 
is a quasi-isomorphism. 
Hence the cofiber of $i$, $\mQ[1] \otimes X$, represents $\Sigma X$.  But
$\mZ[1]$ is invertible with inverse $\mZ[-1]$.  Since the
action of $\ch$ on $\C$ is associative up to coherent isomorphism, this
shows that $\Sigma$ is a self-equivalence on $\Ho(\C)$.  Thus $\C$
is a stable model category.

Since a cylinder object in $\ch$ is also given by tensoring with 
$D^1 \osum \mQ[0]$, 
we have the following natural isomorphisms for a cofibrant object $X$ in 
$\C$ and a fibrant object $Y$ in $\C$.  
\[
H_n\Hom_{\C}(X, Y) \iso [\mZ[n], \Hom_{\C}(X,Y)]^{\Ho(\ch)} \iso 
[\mZ[n] \otimes X, Y]^{\Ho(\C)} \iso 
[\Sigma^n X, Y]^{\Ho(\C)} \] 
\end{proof}

We also need the following definitions from triangulated categories,
see~\cite[1.1]{HPS}.

\begin{definition}
{\em  An object $X$ in a stable model category $\C$ is {\em small} if 
in $\Ho(\C)$ the 
natural map $\Osum_i [X, Y_i] \to [X, \amalg_i Y_i]$ is an
isomorphism for any set of objects $\{Y_i\}$.  A subcategory of a 
triangulated category is {\em localizing} if it is closed under cofiber
sequences, retracts and coproducts. A set of objects $\cG$ in $\C$
is a set of (weak) {\em generators} if the only localizing subcategory of
$\Ho(\C)$ which contains $\cG$ is $\Ho(\C)$.  
}\end{definition}

The following enriched variant of Gabriel's theorem is stated
rationally, but holds integrally as well with $\ch$ replaced by $Ch_{\bZ}$.
A similar statement, on the derived category level, can be found in~\cite{keller}.  Here though we consider the underlying model categories,
before inverting the quasi-isomorphisms, to construct a Quillen
equivalence.

\begin{definition} \label{def-Quillen}
{\em A pair of adjoint functors between model categories is a
{\em Quillen adjoint pair} if the right adjoint preserves fibrations and
trivial fibrations. 
A Quillen adjoint pair induces adjoint total derived functors
between the homotopy categories~\cite[I.4]{Q}.
A Quillen adjoint pair is a {\em Quillen equivalence} if the total
derived functors are adjoint equivalences of the homotopy categories.
This is equivalent to the usual definition by~\cite[1.3.13]{hovey}. }
\end{definition}

\begin{theorem} \label{thm-gabriel} Let $\C$
be a $\ch$-model category and  $\cG$ a set of small, cofibrant and fibrant 
generators. The functor
\[ \HomC (\cG,-) \, : \, {\C} \ \varrow{1cm} \ \Mod\EG \]
is the right adjoint of a Quillen equivalence.
The left adjoint is denoted $-\otimes_{\EG}\cG$.\\
\end{theorem}

\begin{proof}
The functor $\HomC(\cG,-)$ preserves fibrations and trivial fibrations by 
the compatibility axiom (CM) of Definition \ref{def-chain model cat}, since
all objects of $\cG$ are cofibrant.  So together with its left adjoint it 
forms a Quillen pair.  Denote the associated total derived
functors on the homotopy categories by $\RHomC (\cG, -)$ and 
$ - \otimes^L_{\EG}\cG$.  Since $\Ho(\C)$ is a triangulated category both
total derived functors preserve shifts and
triangles; they are exact functors of triangulated categories
by~\cite[I.4 ]{Q}. 
Since $- \otimes^L_{\EG} \cG$ is a left adjoint it commutes with colimits.
The right adjoint also commutes with coproducts since each object of $\cG$
is small.  This follows as in~\cite[3.10.3(ii)]{ss3}.   

We consider the full subcategories of those
$M\in \Ho(\Mod\EG)$ and $X\in\Ho(\C)$ respectively
for which the unit of the adjunction
\[ \eta \ : \ M \ \varrow{1cm} \ \RHomC (\cG, M \otimes^L_{\EG}  \cG) \]
or the counit of the adjunction
\[ \nu \ : \RHomC (\cG,X) \otimes^L_{\EG} X \ \varrow{1cm} \ X \]
are isomorphisms. Since both derived functors are exact and preserve
coproducts, these are localizing subcategories.
For every $G\in\cG$  the $\EG$-module $\HomC (\cG,G)$ is isomorphic to
the free module $F_G^{\EG}$ by inspection
and $F_G^{\EG} \otimes_{\EG} \cG$ is isomorphic to $G$
since they represent the same functor on $\C$.
So the map $\eta$ is an isomorphism for every free module, and the map $\nu$
is an isomorphism for every object of $\cG$.
Since the free modules $F_G^{\EG}$ generate the homotopy category
of $\EG$-modules and the objects of $\cG$ generate $\C$, the localizing
subcategories are the whole categories.  So the  derived
functors are inverse equivalences of the homotopy categories.  
\end{proof}

Since $dg\cA$ is a $\ch$-model category by Proposition~\ref{prop-ch-mc},
we can apply Theorem~\ref{thm-gabriel}.  Triangulated equivalences
preserve generators so we identify generators of $dg\cA$ via the triangulated
equivalence    
$\pi^{\cA} \mc \Ho(\bT\mbox{-spectra}) \to \D(\cA)$ from the homotopy
category of $\bT\mbox{-spectra}$ to the homotopy category of $dg\cA$,
or the derived category of $\cA$, by~\cite[5.6.1]{Greenlees}.
Note that composition of $\pi^{\cA}$ with homology gives the
functor $\pi_*^{\cA}$ which appears in~\cite[5.6.2]{Greenlees}.  The geometric 
basic cells 
$\Bt'=\{\sigma^0_H\}_{H \leq \bT}$ for $H$ any closed subgroup
of $\bT$ form a set of small generators for $\bT$-spectra.  This set
is labeled $\Bt'$ because we will soon consider its image $\Bt$ in
another Quillen equivalent category.  They detect
stable equivalences because $\pi_*^H(X)=[\vee_{K \leq H}\sigma^0, 
X]^{\bT}_*$ for $H$ finite and $\pi_*^{\bT}(X)=[\sigma^0_{\bT}, X]^{\bT}_*$
by~\cite[2.1.5]{Greenlees}.  It follows by~\cite[2.2.1]{ss3} or part
of~\cite[2.3.2]{HPS} that $\Bt'$ is a set of small generators.
Hence, the set of 
images $ \{\pi^{\cA}(\sigma^0_H) = L_H \}_{H \leq \bT}$ form a set of 
small generators for $dg\cA$.   Let $\Ba= \{ I_H \}_{H \leq \bT}$ be
the set of their fibrant replacements.  Then $\Ba$ is a set of small,
cofibrant and fibrant generators.  Let $\E(\Ba)$ be the endomorphism
$\ch$-category as in Definition~\ref{def-EP}.  
In Section~\ref{sec-endo} we explicitly describe $\E(\Ba)$. 
 
\begin{corollary}\label{cor-Gabriel}
The model category $dg\cA$ is Quillen equivalent to the model category
of differential graded modules over $\E(\Ba)$.
\end{corollary}

To prepare for constructing a Quillen equivalence between $dg\cA$
and the model category of rational $\bT$-equivariant spectra, $\bT$-spectra,
we show that there is an equivalence between $\bT$-spectra and
modules over a rational differential graded ring with many objects.
A Morita equivalence similar to Theorem~\ref{thm-gabriel} in~\cite{ss3} 
considers $\Spec^{\Sigma}$-model categories, or 
{\em spectral model categories}, 
model categories compatibly enriched over symmetric spectra.
\cite[3.10.3]{ss3} shows that any spectral model category
with a set of small generators
is equivalent to modules over a ring spectrum with many objects.  In this paper
we are interested in particular in {\em rational stable model categories}, 
stable model categories where $[X,Y]^{\Ho(\C)}$ is a rational vector space for
any objects $X$ and $Y$ in $\C$.  In this case the Gabriel equivalence
produces a rational ring spectrum which can be replaced by a differential
graded ring with many objects.

\begin{theorem}~\cite[1.3]{ssme}\label{thm-rational-Gabriel}
Let $\C$ be a spectral model
category.  If $\C$ is rational and has a set of small generators, then
$\C$ is Quillen equivalent (via a chain of several Quillen equivalences) to 
the model category of modules over a rational 
differential graded ring with many objects.
\end{theorem}

Before rationalization, the category of $\bT$-equivariant spectra  is
a cellular, simplicial, left proper, 
stable model category ~\cite[III.4.2]{MM}.  In particular,
a $G$-topological model category~\cite[III.1.14]{MM} is a topological 
model category which in turn is a simplicial model category via the 
singular and realization functors.  Also, a compactly generated model
category is a stricter notion than a cofibrantly generated model category.
One can check that $\bT$-equivariant spectra is in fact cellular, see
~\cite[A.1]{hovey-spec} or~\cite[14.1.1]{psh}.  Since localizations
preserve these properties by~\cite[4.1.1]{psh}, the model category of rational 
$\bT$-equivariant spectra is also a cellular, simplicial, left proper,
stable model category.  Then Hovey in~\cite[8.11,9.1]{hovey-spec} shows
that $\bT$-spectra is Quillen equivalent to a spectral model category, $\C'$.
Since $\bT$-spectra is also a simplicial, cofibrantly generated, proper,
stable model category,~\cite[3.9.2]{ss3} also shows it is
Quillen equivalent to a spectral model category. 
Since rational $\bT$-equivariant spectra has a set of small generators, 
$\B_t'$, the geometric basic cells, $\C'$ is also
rational and has a set of small generators because these properties are
determined on the homotopy category level.
Thus, by Theorem~\ref{thm-rational-Gabriel}, $\C'$ is Quillen
equivalent to $\Mod\cR$ for some rational differential graded
ring with many objects, $\cR$.  Let $\BCt$ be a set of cofibrant
and fibrant replacements of the images of the basic cells 
in $\Mod\cR$.  Then $\BCt$ is a set of small generators in $\Mod
\cR$. 
Since $\Mod\cR$ is a $\ch$-model category, Theorem~\ref{thm-gabriel} shows that 
$\Mod\cR$ is Quillen equivalent to $\Mod\E(\BCt)$ where $\E(\BCt)$ is the 
endomorphism ring on the generators $\BCt$.

\begin{corollary}\label{cor-Gabriel-t}
$\bT$-spectra is Quillen equivalent to $\Mod\E(\BCt)$ which is the 
endomorphism ring on a set of small generators $\BCt$ in an intermediate
model category $\Mod\cR$.
$$ \bT\mbox{-spectra} \Qeq \Mod\cR \Qeq \Mod\E(\BCt)$$
\end{corollary}

This does give an algebraic model for $\bT$-spectra, but $\E(\BCt)$ is 
large and not explicit.   This algebraic model is nonetheless useful
because we can construct a zig-zag of quasi-isomorphisms between $\E(\BCt)$
and $\E(\Ba)$ which induce Quillen equivalences on the respective module 
categories.  

\section{The endomorphism ring for $dg\cA$}\label{sec-endo}

In this section we first explicitly calculate $\E(\Ba)$, the endomorphism
ring associated to $\dga$ by Corollary~\ref{cor-Gabriel}.  Then in
Definition~\ref{def-Ea} we define a sub-ring $\Ea$ which is quasi-isomorphic 
to $\E(\Ba)$.  In Section~\ref{sec-map}, we then use this simpler ring
$\Ea$ to construct
a zig-zag of Quillen equivalences between $\dga$ and $\bT$-spectra.  
For this construction, we need to understand the Massey products in 
$H_*\Ea$, so we end this section with a definition of Massey products
for differential graded rings with many objects and a calculation of the 
products and Massey products in $H_*\Ea$.  Since $\Ea$ and $\E(\Ba)$ are
quasi-isomorphic this calculation also applies to $H_*\E(\Ba)$ 
by~\cite[1.5]{may}.

First we calculate
the underlying rational chain complexes of $\E(\Ba)$.  

\begin{proposition}\label{prop-vector space}
For any distinct finite subgroups $H, K$ of $\bT$
the rational chain complexes of homomorphisms are as follows:
\begin{enumerate}
\item $\Hom_{\dga}(I_H, I_H) = (\Osum_{n \leq 0}D^n) \osum \bQ[0] \osum \bQ[1]$,  
\item $\Hom_{\dga}(I_H, I_K) = 0$,
\item $\Hom_{\dga}(I_H, I_{\bT}) = (\Osum_{n \leq 0} D^{2n -1}) \osum \bQ[0]$, 
\item $\Hom_{\dga}(I_{\bT}, I_H) = (\Osum_{n \geq 0} D^{2n+1}) \osum
(\Osum_{m \leq 0} D^m) \osum \bQ[1]$ and 
\item $\Hom_{\dga}(I_{\bT}, I_{\bT}) = (\Osum_{n \geq 0} \bQ\cF[2n +1]) 
\osum (\Osum_{m \leq 0} (D^{2m}\otimes \bQ\cF)) \osum \bQ[0]$ where 
$\bQ\cF$
is the rational vector space with basis $\cF$, the set of finite subgroups
of $\bT$.
\end{enumerate}
\end{proposition} 

\begin{proof}
In the first three cases, since the vertex of the source is trivial the maps 
are determined by considering graded $\cO_{\cF}$-module maps of the nubs,
denoted by $\GrO(-,-)$.  
As an underlying $\cO_{\cF}$-module, the nub of $I_H$ is 
$\Sigma \bI(H) \osum \Sigma^2
\bI(H)$.  As a graded vector space $\GrO(\bI(H), \bI(H)) = 
\Osum_{n \leq 0} \bQ[2n]$.
This determines the underlying graded vector space of $\Hom_{\dga}(I_H, I_H)$.
Let $\{m_{2n}\}_{n \leq 0}$ be a basis for $\GrO(\Sigma\bI(H), \Sigma\bI(H))$,
$\{m_{2n +1}\}_{n \leq 0}$ be a basis for $\GrO(\Sigma\bI(H), \Sigma^2\bI(H))$,
$\{l_{2n}\}_{n \leq 0}$ be a basis for 
$\GrO(\Sigma^2\bI(H), \Sigma^2\bI(H))$ and
$\{l_{2n-1}\}_{n \leq 0}$ be a basis for $\GrO(\Sigma^2\bI(H), \Sigma\bI(H))$.
When necessary, these bases and those below are decorated with an $H$.   
The differential
is determined by $dm_{2n} = m_{2n -1}$, $dm_{2n+1} = 0$, $dl_{2n}=m_{2n -1}$
and $dl_{2n-1}=l_{2n-2} - m_{2n-2}$.  
Since $\GrO(\bI(H), \bI(K))=0$, it follows that
$\Hom_{\dga}(I_H, I_K)= 0.$  
 
As an underlying $\cO_{\cF}$-module, the nub of $I_{\bT}$ is 
$ t_*^{\cF} \osum \Sigma \bI$.
As a graded vector space $\GrO(\bI(H), \bI)=
\Osum_{n \leq 0} \bQ[2n]$ and $\GrO(\bI(H), t_*^{\cF})=0$;
this determines the underlying graded vector space of $\Hom_{\dga}(I_H, 
I_{\bT})$.  Let $\{f_{2n}\}_{n \leq 0}$ be a basis for 
$\GrO(\Sigma\bI(H), \Sigma\bI)$ and $\{f_{2n -1}\}_{n \leq 0}$ be a basis for 
$\GrO(\Sigma^2\bI(H), \Sigma\bI)$.  The differential is
determined by $df_{2n-1}=-f_{2n-2}$ and $df_{2n}=0$.

For the fourth case the vertex of the target is trivial, so again
the maps are determined by the nubs.  Also, note that 
$\GrO(t_*^{\cF}, \bI(H))= \Osum_{n \in \bZ} \bQ[2n]$.     
Let $\{g_{2n+1}\}_{n \in \bZ}$ be a basis for $\GrO(t_*^{\cF}, \Sigma\bI(H))$,
$\{g_{2n }\}_{n \in \bZ}$ be a basis for $\GrO(t_*^{\cF}, \Sigma^2\bI(H))$,
$\{h_{2n}\}_{n \leq 0}$ be a basis for 
$\GrO(\Sigma\bI(H), \Sigma\bI(H))$ and
$\{h_{2n+1}\}_{n \leq 0}$ be a basis for 
$\GrO(\Sigma\bI(H), \Sigma^2\bI(H))$.   The differential is
determined by $dg_{2n+1}= g_{2k}$, $dg_{2n} = 0$, $dh_{2n}= h_{2n-1} + 
g_{2n-1}$ and $dh_{2n +1}= -g_{2n}$. 

The last case is the only one where we must consider the vertex.  This
restricts the maps from $t_*^{\cF} \otimes \bQ$ to itself to maps of
the form $1 \otimes \phi$.  Hence $\GrO(t_*^{\cF}, t_*^{\cF})$ 
contributes only $\bQ[0]$
to $\Hom_{\dga}(I_{\bT}, I_{\bT})$.  Also, note that $\GrO(t_*^{\cF}, \bI)=
\Osum_{n \in \bZ}\bQ\cF[2n]$ and $\GrO(\bI, \bI)= \Osum_{n \leq 0}
\bQ\cF[2n]$.  Let $\{i^H_{2n+1}\}_{n \in \bZ, H \in \cF}$ be a basis for 
$\GrO(t_*^{\cF}, \Sigma\bI(H))$,
$\{j^H_{2n}\}_{n \leq 0, H \in \cF}$ be a basis for 
$\GrO(\Sigma\bI(H),\Sigma\bI(H))$  and $\id_{\bT}$ the identity 
map.   The differential is determined by
$di^H_{2n+1}=0$, $dj^H_{2n}= i^H_{2n-1}$ and $d\id_{\bT} =0$.
\end{proof}

As we note below in Proposition~\ref{prop-positive}, for the rest of this
paper we only need to understand the
product structure on the $(-1)$-connected cover of $\E(\Ba)$.  This
greatly simplifies the structure, see Proposition~\ref{prop-one}.    
So only the non-negative dimensional products in the following proposition
are needed and these products are displayed in Definition~\ref{def-Ea}.

\begin{proposition}\label{prop-products}
The following is a list of all of the non-trivial products, or compositions, 
of the bases chosen for $\E(\Ba)$ in Proposition~\ref{prop-vector space}.
\[\begin{array}{llll}
m_{2n} g_{2k +1}= g_{2n +2k +1}&  
m_{2n+1} g_{2k +1}= g_{2n +2k +2}&  
l_{2n} g_{2k} = g_{2n +2k} & 
l_{2n-1} g_{2k} = g_{2n +2k -1}\\
f_{2n-1} g_{2k} = i_{2n +2k -1}&  
f_{2n} g_{2k+1} = i_{2n +2k +1}&
l_{2n} h_{2k} = g_{2n +2k } & 
l_{2n-1} h_{2k} = g_{2n +2k-1 }\\  
l_{2n} h_{2k+1} = h_{2n +2k+1 }&  
l_{2n-1} h_{2k+1} = h_{2n +2k }&  
h_{2n} f_{2k} = m_{2n +2k } & 
h_{2n} f_{2k-1} = l_{2n +2k-1 }\\  
h_{2n+1} f_{2k} = m_{2n +2k+1 }&  
h_{2n+1} f_{2k-1} = l_{2n +2k }& 
f_{2n} h_{2k} = j_{2n +2k } & 
f_{2n-1} h_{2k+1} = j_{2n +2k }\\  
h_{2n} j_{2k} = h_{2n +2k } & 
h_{2n+1} j_{2k} = h_{2n +2k+1 }& 
h_{2n} i_{2k+1} = g_{2n +2k+1 }&  
h_{2n+1} i_{2k+1} = g_{2n +2k+2 }\\  
f_{2n} m_{2k} = f_{2n +2k } & 
f_{2n-1} m_{2k+1} = f_{2n +2k }&  
f_{2n-1} l_{2k} = f_{2n +2k-1 }&  
f_{2n} l_{2k-1} = f_{2n +2k-1 }\\  
j_{2n} f_{2k} = f_{2n +2k } & 
j_{2n} f_{2k-1} = f_{2n +2k-1 }&  
m_{2n} m_{2k} = m_{2n +2k } & 
m_{2n+1} m_{2k} = m_{2n +2k +1}\\  
m_{2n} l_{2k-1} = l_{2n +2k -1}&  
m_{2n+1} l_{2k-1} = l_{2n +2k }& 
l_{2n} m_{2k+1} = m_{2n +2k+1 }&  
l_{2n-1} m_{2k+1} = m_{2n +2k }\\  
l_{2n} l_{2k} = l_{2n +2k } & 
l_{2n-1} l_{2k} = l_{2n +2k-1 }&  
j^H_{2n} j^H_{2k} = j^H_{2n +2k }&  
j^H_{2n} i^H_{2k+1} = i^H_{2n +2k+1 }\\ 
\end{array}
\]
Also, $x \id_{\bT} = x$ {where} $x=i_{2n+1}, j_{2n}, \id_{\bT}$ or $f_n$
and $\id_{\bT} y= y$ {where} $y=i_{2n+1}, j_{2n}, \id_{\bT}, g_n$ or $h_n$.
\end{proposition}

Next we consider simpler sub-rings of $\E(\Ba)$.  

\begin{definition}\label{def-positive}
{\em The {\em $(-1)$-connected cover} $R\langle 0 \rangle$ of a differential 
graded ring $R$ is defined by $R\langle 0 \rangle_n= 0$ for $n < 0$,
$R\langle 0 \rangle_n = R_n$ for $n > 0$ and $R\langle 0 \rangle_0
=Z_0 R$ where $Z_0 R$ is the kernel of the differential $d\mc R_0 \to R_{-1}$.
Since the product of two cycles is a cycle, the product on $R$ induces
a product on $R\langle 0\rangle$.  Moreover the inclusion $R\langle 0 \rangle
\to R$ is a map of differential graded rings which induces an isomorphism
in homology in non-negative degrees.  This definition generalizes to a 
differential graded ring with many objects, $\cR$, by applying this to each morphism
complex $\cR(o,o')$.  
}\end{definition}

Since the homology of $\E(\Ba)$ is concentrated in non-negative degrees
it is quasi-isomorphic to its $(-1)$-connected cover.  By~\cite[4.3]{ss1},
the module categories are Quillen equivalent. 

\begin{proposition}\label{prop-positive}
The map $\E(\Ba)\langle 0 \rangle \to \E(\BC)$ is a quasi-isomorphism
of differential graded rings with objects $\BC$.
Thus, there is a Quillen equivalence between the respective categories
of modules.
\[ \Mod\E(\BC)\langle 0\rangle  \Qeq \Mod\E(\BC)\]
\end{proposition}

Next we explicitly describe $\E(\BC)\langle 0 \rangle$.
\begin{proposition}\label{prop-one}
The rational chain complexes underlying the differential graded ring 
$\E(\BC)\langle 0\rangle$ are given by: 
\begin{enumerate}
\item $\EBCo(I_H, I_H) = \bQ[0] \osum \bQ[1]$ with generators $\id_H = m_0 + l_0$
and $m_1$, 
\item $\EBCo(I_H,I_K) = 0$,
\item $\EBCo(I_H, I_{\bT}) = \bQ[0]$ with generator $f_0$, 
\item $\EBCo(I_{\bT}, I_H) = (\Osum_{n \geq 0} D^{2n+1}) \osum \bQ[1]$ with
generators $g_n$ for degrees $n \geq 2$, $g_1 + h_1$ generating the
cycles in degree $1$ and  $h_1$ and $g_0$ generating $D^1$ and 
\item $\EBCo(I_{\bT}, I_{\bT}) = (\Osum_{n \geq 0} \bQ\cF[2n +1])
\osum \bQ[0]$ with generators $i^H_{2n+1}$ and $\id_{\bT}$ with $n \geq 0$. 
\end{enumerate}
The products are determined by the products in $\E(\BC)$.
\end{proposition} 

We make one more modification to the ring $\E(\BC)\langle 0 \rangle$ 
by considering the sub-ring $\Ea$. 

\begin{definition}\label{def-Ea}
{\em 
Define the differential graded ring $\Ea$ with objects $\BC$ as the sub-ring of 
$\E(\BC)\langle 0 \rangle$ produced by deleting the $D^1$ generated by $h_1$ 
and $g_0$.  Denote the element $g_1 + h_1$
in $\E(\BC)\langle 0 \rangle$ by $\tg_1$ in $\Ea$.
To ease notation in Section~\ref{sec-map}, the element $\tg_1$ will
be renamed $g_1$ to match the family of elements $g_{2n+1}$.
Then the rational chain complexes underlying $\Ea$ are given by: 
\begin{enumerate}
\item $\Ea (I_H, I_H) = \bQ[0] \osum \bQ[1]$ with generators $\id_H = m_0^H + l_0^H$
and $m_1^H$, 
\item $\Ea (I_H, I_K) = 0$,
\item $\Ea (I_H, I_{\bT}) = \bQ[0]$ with generator $f_0^H$, 
\item $\Ea (I_{\bT}, I_H) = (\Osum_{n \geq 1} D^{2n+1}) \osum \bQ[1]$ with
generators $\tg_1^H$, $g_n^H$ for $n \geq 2$, and
\item $\Ea (I_{\bT}, I_{\bT}) = (\Osum_{n \geq 0} \bQ\cF[2n +1])
\osum \bQ[0]$ with generators $i_{2n+1}^H$ and $\id_{\bT}$ with $n \geq 0$. 
\end{enumerate}
The non-trivial products of these chosen bases, except the obvious products with
$\id_H$ and $\id_{\bT}$ are:
\begin{enumerate}
\item $\tg_1^H f_0^H  = m_1^H$, 
\item $f_0^H g_{2n+1}^H = i_{2n+1}^H$ for $n \geq 0$ and 
\item $\tg_1^H f_0^H g_{2n+1}^H = g_{2n+2}^H$ for $n \geq 0$. 
\end{enumerate}
When $n=0$, the second line should be read as $f_0 \tg_1^H = i_1^H$ and
the third line should be read as $\tg_1^H f_0^H \tg_1^H = g_2^H$.
The products in the third line can also be rewritten in terms of $m_1^H$ or 
$i_{2n+1}^H$ using the first two lines.  
}\end{definition} 

Another description of $\Ea$ is that it is the differential graded
ring generated by the elements $f_0^H$, $\tg_1^H$ and $g_{2n+1}^H$ for $n > 0$
with relations generated by $f_0^H \tg_1^H f_0^H = 0$, 
$g_{2n+1}^H f_0^H = 0$ for $n > 0$ and 
$x^H y^K =0$ whenever $H \neq K$. 
These relations are equivalent to requiring $\Ea(I_H, I_{\bT}) \iso \mQ[0]$,
$\Ea(I_H, I_H) \iso \mQ[0] \osum \mQ[1]$ and $\Ea(I_H, I_K) =0$.

Since quasi-isomorphisms induce Quillen equivalences 
by~\cite[4.3]{ss1}, we have the following statement.   

\begin{proposition}\label{prop-simple-alg}
The inclusion of differential graded rings 
$\Ea \to \E(\BC)\langle 0 \rangle$ which sends $\tg_1$ to $g_1 + h_1$ is a 
quasi-isomorphism.  Hence the composite $\Ea \to \EBCo \to \E(\BC)$ is
also a quasi-isomorphism.
Thus, there is a Quillen equivalence between the respective categories
of modules.
\[ \Mod\Ea \Qeq  \Mod \E(\BC)\]
\end{proposition}

We now consider higher order products. Here the higher order products
we consider are Massey products. 
We generalize the definition
of Massey products for a differential graded algebra due to 
Massey~\cite{massey} to 
a differential graded ring with many objects; see also~\cite[8.3.2]{mccleary}.
Since there is no
real difference between the homological and cohomological versions we
still refer to these as Massey products instead of the `Yessam' products
of Stasheff~\cite{stasheff}.  
Notice that with
this definition the indeterminacy may be smaller for a differential
graded ring with many objects than for the associated differential
graded matrix ring on only one object.  To deal with signs here
we define  $\epsilon(a)=(-1)^{\mbox{\scriptsize deg } a + 1}a$ instead
of using the usual bar notation.

\begin{definition}\label{def-Massey}
{\em Let $\cR$ be a differential graded ring with objects 
$\cG=\{G_i\}_{i \in I}$.  
Suppose that $\gamma_1, \dots, \gamma_n$ are classes in $H_*\cR$ with
$\gamma_i \in H_{p_i}\cR(G_i, G_{i-1})$.  A {\em defining system},
associated to $\langle \gamma_1, \dots, \gamma_n\rangle$, is
a set of elements $a_{i,j}$ for $1\leq i\leq j \leq n$ with $(i,j) \neq (1,n)$
such that 
\begin{enumerate}
\item $a_{i,j} \in \cR(G_j, G_{i-1})_{p_i + \dots + p_j + j -i}$,
\item $a_{i,i}$ is a cycle representative of $\gamma_i$ and 
\item $d(a_{i,j})= \Sigma_{r=i}^{j-1}  \epsilon(a_{i,r}) a_{r+1, j}$ with
$\epsilon(a)=(-1)^{\mbox{\scriptsize deg } a + 1}a$.
\end{enumerate}
To each defining system associate the cycle
\[ \Sigma_{r=1}^{n-1} \epsilon(a_{1,r}) a_{r+1,n} \in \cR(G_n, G_0)_{p_1 + \dots
p_n + n - 2} \]  
The {\em $n$-fold Massey product}, 
$\langle \gamma_1, \dots, \gamma_n \rangle$, is the set of all homology
classes of cycles associated to all possible defining systems. 
}\end{definition}

Since Massey products are preserved
by quasi-isomorphisms by~\cite[1.5]{may}, 
we can use the simpler ring $\Ea$ to calculate these products for $\E(\BC)$.

\begin{proposition}\label{prop-homology-brackets}
The homology of $\Ea$ (and hence also of $\E(\BC)$ and $\EBCo$) is given
by
\begin{enumerate}
\item $H_*\Ea (I_H, I_H) = \bQ[0] \osum \bQ[1]$ with generators $[\id_H]$
and $[m_1^H]$, 
\item $H_*\Ea (I_H, I_K) = 0$,
\item $H_*\Ea (I_H, I_{\bT}) = \bQ[0]$ with generator $[f_0^H]$, 
\item $H_*\Ea (I_{\bT}, I_H) = \bQ[1]$ with
generator $[\tg_1^H]$, and 
\item $H_*\Ea (I_{\bT}, I_{\bT}) = (\Osum_{n \geq 0} \bQ\cF[2n +1])
\osum \bQ[0]$ with generators $[i^H_{2n+1}]$ and $[\id_{\bT}]$. 
\end{enumerate}
The non-trivial products, except the obvious ones with $[\id_H]$ and 
$[\id_{\bT}]$, and the non-trivial Massey products are given by
\begin{enumerate}
\item $[\tg_1^H] [f_0^H] = [m_1^H]$,
\item $[f_0^H] [\tg_1^H] = [i_1^H]$, and
\item  $\langle [f_0^H], \underbrace{[m_1^H], \dots, [m_1^H]}_n, [\tg_1^H] \rangle = 
\{[-i^H_{2n+1}]\}$.
\end{enumerate}
\end{proposition}
\noindent Note that $[\tg_1^H]$ and $[f_0^H]$ generate all of the elements in 
homology via products and Massey products.

\begin{proof}
The products follow from the fact that $\tg_1 f_0 = m_1$ and $f_0 \tg_1 = i_1$.
Next consider the triple product $\langle [f_0^H], [m_1^H], [\tg_1^H] \rangle$.
The cycle representatives here are unique. 
Since $f^H_0 m_1^H = 0$, $\Ea (I_H,I_\bT)_2 =0$ and 
$g_3^H$ 
is the unique element such that $dg^H_3=g_2^H= m_1^H \tg_1^H$, this 
determines the only defining system as the following matrix.
\[
\left(\begin{array}{ccc}
f_0^H & 0 & \\
&m_1^H & g_3^H \\
&&\tg_1^H
\end{array}
\right)
\] 
The associated cycle is $-f^H_0 g^H_3 +0 \tg_1^H = -i^H_3$, so
$\langle [f_0^H], [m_1^H], [\tg_1^H] \rangle = \{[-i^H_3]\} $.

The higher order products follow inductively.  At each stage because 
$\Ea(I_H, I_H)_n$ and $\Ea(I_H,I_{\bT})_n$ are trivial for $n > 1$ and
$\Ea(I_{\bT}, I_H)_n$ is generated by one element, 
there are no choices and hence no indeterminacy.
Set $a_{1,1}= 
f_0^H$, $a_{i,i}= m_1^H$ for $ 2 \leq i \leq n+1$, and $a_{n+2,n+2}= \tg_1^H$.
Since $f^H_0 m_1^H = 0$ and $m_1^H m_1^H =0$, the only choice for   
$a_{i,j}$ is $0$ for $i < j$ and $j \neq n+2$.  
Hence the only non-zero entries in the defining system are $a_{i,i}$ 
for $1 \leq i \leq n+2$ and 
$a_{i,n+2}$ for $1 < i \leq n+2$.  Then the only choice
for $a_{i, n+2}$ is $g_{2n-2i+5}^H$ for $i < n+2 $.  
Since $-a_{1,1} a_{2,n+2}=-f_0^H g^H_{2n+1}
=-i^H_{2n+1}$, the conclusion follows.
\end{proof}

\section{The endomorphism ring for $\bT$-spectra}\label{sec-end-t}

In this section we consider the endomorphism ring $\EBCt$ associated
to $\bT$-spectra in Corollary~\ref{cor-Gabriel-t}.  As with $\E(\Ba)$,
we modify $\EBCt$ to define a quasi-isomorphic ring $\Et$. 
One property of $\Et$ is that $\Et(H,\bT)$ is concentrated
in degree zero.  This forces certain relations to hold and thus makes
it easier to construct a zig-zag of quasi-isomorphisms between $\Ea$
and $\Et$.  Another property of $\Et$ is that in degree zero
$(\Et)_0 \iso H_0\Et$. 
This gives a place to start the construction of
the maps. This is carried out in Section~\ref{sec-map}. 

We first modify $\EBCt$ by taking its $(-1)$-connected cover
$\EBCto$ as in Definition~\ref{def-positive}.  We next modify $\EBCto$ so
that degree zero agrees with its homology.  This uses 
a Postnikov approximation (Definition~\ref{def-postnikov}); its
general properties are stated in Propositions~\ref{prop-5.4} and~\ref{prop-bar}
and are applied to this specific case in Proposition~\ref{prop-pullback}.
Finally, we define $\Et$ as a sub-ring of the resulting differential
graded ring with many objects.  This ring $\Et$ then has the properties 
mentioned above.

We first recall the homology of $\EBCt$ from~\cite[2.1.6]{Greenlees}.  Of
course, it is isomorphic as a ring to the homology of 
$\EBCa$ stated in Proposition
\ref{prop-homology-brackets} because the triangulated equivalence 
of $\Ho(\bT\mbox{-spectra})$ and $\cD(\cA)$
from \cite[5.6.1]{Greenlees} can be chosen in such a way that the 
geometric basic cells in $\Bt$ are taken to the algebraic basic cells in
$\Ba$.  Because of this isomorphism we use the notation $\wh{x}$ to denote
the element in $H_*\E(\Bt)$ corresponding to $[x]$ in $H_*\Ea$.  We also
delete most of the superscripts $H$.  We label the elements in $\Bt$ by their 
corresponding subgroup in $\bT$.

\begin{proposition}\label{prop-homology-top}
The homology of $\EBCt$ is given by 
\begin{enumerate}
\item $H_*\EBCt(H, H) = \bQ[0] \osum \bQ[1]$ with generators $\wh{\id_H}$ and
$\wh{m_1}$,  
\item $H_*\EBCt(H, K) = 0$,
\item $H_*\EBCt(H, {\bT}) = \bQ[0]$ with generator $\wh{f_0}$, 
\item $H_*\EBCt({\bT}, H) = \bQ[1]$ with generator $\wh{g_1}$ and 
\item $H_*\EBCt({\bT}, {\bT}) = (\Osum_{n \geq 0} \bQ\cF[2n +1])
\osum \bQ[0]$ with generators $\wh{i_{2n+1}^H}$ and $\wh{\id_{\bT}}$. 
\end{enumerate}
The non-trivial products, except the obvious ones with $\wh{\id_H}$ and 
$\wh{\id_{\bT}}$ are given by
\begin{enumerate}
\item $\wh{g_1^H} \wh{f_0^H} = \wh{m_1^H}$,
\item $ \wh{f_0^H} \wh{g_1^H} = \wh{i_1^H}$, and
\end{enumerate}
\end{proposition}

The next step is to consider the $(-1)$-connected cover of $\EBCt$.
Since the homology of $\EBCt$ is concentrated in non-negative degrees,
$\EBCto \to \EBCt$ is a quasi-isomorphism.
\begin{proposition}\label{prop-positive-top} 
The map $\EBCto  \to \EBCt$ is a quasi-isomorphism
of differential graded rings with objects $\Bt$.
Thus, there is a Quillen equivalence between the respective categories
of modules.
\[ \Mod\EBCto  \Qeq \Mod\EBCt\]
\end{proposition}

The purpose of the next step is to construct a ring $\overline{\cR}$
quasi-isomorphic to $\EBCto$ such that $\overline{\cR}_0 \iso 
H_0(\overline{\cR})$.  This involves a Postnikov
approximation which is obtained by killing homology
groups above degree zero.  The next definition gives a functorial
definition of this approximation which uses the small object argument; 
see~\cite[7.12]{dwyer-spalinski}.

Consider the functor $\ev_{(H,K)}$ from rational differential graded rings
with objects $\B_t$ to rational chain complexes which is evaluation
at $(H,K)$.  Denote the left adjoint of $\ev_{(H,K)}$ by $T_{(H,K)}$.
$T_{(H,H)}M$ is the tensor algebra on $M$ at $(H,H)$, $\bQ[0]$ at $(K,K)$
for $H \neq K$ and trivial everywhere else.  
 
\begin{definition}\label{def-postnikov}
{\em Given a rational differential graded ring $\cR$ with objects 
$\B_t$ a  
Postnikov approximation of $\cR$ is given by the factorization $\cR \varr{i} 
P_0\cR \to 0$
produced by the small object argument for the set of maps $I=\{ T_{(H,K)}
(\bQ[n] \to D^{n+1})\}$ for $n \geq 1$ and all pairs $(H,K) \in \B_t \times 
\B_t$.
}\end{definition}

Since the map $i\mc \cR \to P_0\cR$ is constructed by pushouts of the maps
in $I$ and $P_0\cR \to 0$ has the right lifting property with respect to
the maps in $I$, these maps have the following properties.

\begin{proposition}\label{prop-5.4}
The map $\cR \varr{i} P_0\cR$ is an isomorphism in degrees less than or equal 
to one and hence induces an isomorphism in homology in non-positive degrees.  
The map $P_0\cR \to 0$ induces an isomorphism in homology in
positive degrees.
\end{proposition}

By applying the small object argument one more time we factor the
map $\cR \varr{i} P_0\cR$.  Construct the factorization $\cR 
\varr{\bar{i}} \overline{\cR}
\varr{\bar{p}}P_0\cR$ by using the small object argument with respect to 
the set $J= \{T_{(H,K)}(0 \to D^n)\}$ for $n \geq 2$ and all pairs
$(H,K) \in \B_t \times \B_t$.  Again by construction, since the maps
in $J$ are injective quasi-isomorphisms and isomorphisms in non-positive
degrees and the map $\cR \varr{i} P_0\cR$ has the above listed properties, 
$\bar{i}$ and $\bar{p}$ have the following properties.

\begin{proposition}\label{prop-bar}
The map $\cR\varr{\bar{i}} \overline{\cR}$ is a quasi-isomorphism which is an 
isomorphism in non-positive degrees.  $\overline{\cR} \varr{\bar{p}} P_0\cR$
is an epimorphism and in non-positive degrees is a quasi-isomorphism and an 
isomorphism. 
\end{proposition}

\begin{proof}
Since $\bar{p}$ has the lifting property with respect to $J$, 
$\ev_{(H,K)}\bar{p}$ is an epimorphism in degrees greater than one.  
Since $\cR \varr{i} P_0\cR$
is an isomorphism in degrees less than two and $i=\bar{p} \bar{i}$, it follows
that $\bar{p}$ is also an epimorphism in degrees less than two.
\end{proof}

We now apply these general constructions to $\EBCto$.  For ease of notation,
let $\E_t'$ stand for $\EBCto$. Since $\E_t'$ is
$(-1)$-connected, it is trivial in negative degrees and its homology
is concentrated in non-negative degrees.   Hence, $P_0\E_t'$ 
is trivial in negative degrees and its homology is concentrated in degree zero. 
Let $\E_t' \varr{\bar{i}} \overline{\E_t'} \varr{\bar{p}}
P_0\E_t'$ be the factorization constructed above for $\cR=\E_t'$. 
Consider $H_0\E_t'\iso H_0\EBCt$ as a differential graded ring with 
objects $\Bt$ which is concentrated in degree zero.  
From Proposition~\ref{prop-homology-top}, $H_0\EBCt$ is the free $\mQ$-algebra
with objects $\B_t$ generated by 
a one dimensional vector space of morphisms from $H$ to $\bT$ for each
$H \in \cF$.   
Since $H_0\EBCt \iso H_0\Et' \iso H_0P_0\E_t'$, one can construct 
a quasi-isomorphism of differential graded rings $q\mc H_0\E_t' \to P_0\E_t'$.  
Define $\E_t''$ as the differential graded ring with objects
$\Bt$ which is the pullback of $\bar{p}$ and $q$.
\[\begin{CD}
\E_t'' @>{q'}>> \overline{\E_t'}\\
@VVV @VV{\bar{p}}V\\
H_0\E_t' @>>{q}> P_0\E_t'
\end{CD}\]

\begin{proposition}\label{prop-pullback}
The map $q'\mc \E_t'' \to \overline{\E_t'}$ is a quasi-isomorphism.  Moreover,
$(\E_t'')_0 \iso H_0\E_t'' \iso H_0\E_t'$ and $(\E_t'')_n= 0$ for
$n < 0$.  
\end{proposition}

\begin{proof}
$\E_t''(H,K)$ is the pullback of differential graded modules
of the maps $\ev_{(H,K)}\bar{p}$ and $\ev_{(H,K)}q$. 
Each map $\ev_{(H,K)}{\bar{p}}$ is an epimorphism by Proposition~\ref{prop-bar}
and each map $\ev_{(H,K)}q$ is a quasi-isomorphism.  
The model category of differential graded modules is right proper since every 
object is fibrant.  So $\E_t''(H,K) \to \overline{\E_t'}(H,K)$ is a 
quasi-isomorphism 
because it is the pullback of a weak equivalence across a fibration. 
Since ${\bar{p}}$ is an isomorphism in non-positive degrees by 
Proposition~\ref{prop-bar}, $(\E_t'')_n$ is isomorphic to $(H_0\E_t')_n$
for $n \leq 0$.  Note that since $q'$ is a quasi-isomorphism this
also implies that $d_1 \mc (\E_t'')_1 \to (\E_t'')_0$ is trivial.
\end{proof}

We make one more modification to $\E_t''$.  Define $\Et$ as the
sub-differential graded ring with objects $\Bt$ which differs from
$\E_t''$ only by setting $\Et(H, \bT) = \E_t''(H,\bT)_0$ concentrated in 
degree zero for each $H \in \cF$ and $\Et(H,K) =0$ for $H,K$ distinct
finite subgroups.
This defines a differential graded ring since $\E_t''(H,\bT)_0$
is a retract of $\E_t''(H,\bT)$ because $d_1$ is trivial. 
It follows that the map $\Et \to \E_t''$ is a quasi-isomorphism and
$\Et$ has the properties listed in Proposition~\ref{prop-pullback}
for $\E_t''$.  Combining this with Propositions~\ref{prop-positive-top} 
and~\ref{prop-bar} gives the following.

\begin{proposition}\label{prop-simple-top}
The composite map $\Et \to \E_t'' \to \overline{\E_t'}$
is a quasi-isomorphism of differential graded rings with objects $\Bt$.
Combining this with the quasi-isomorphisms $\EBCto \iso \E_t' \to 
\overline{\E_t'}$ and $\EBCto \to \EBCt$ gives a zig-zag of quasi-isomorphisms
between $\Et$ and $\EBCt$.
Thus, there is a zig-zag of Quillen equivalences between the respective 
categories of modules.
\[ \Mod\Et \Qeq \Mod\EBCt\]
Moreover, $(\Et)_0\iso H_0(\Et)$, $(\Et)_n=0$ for $n < 0$, $\Et(H,K) =0$ and
$\Et(H,\bT)$ is concentrated in degree zero.
\end{proposition}

\section{The quasi-isomorphism between $\Ea$ and $\Et$}\label{sec-map}

In this section we construct a differential graded ring $\cS$ with 
objects $\Ba$ and a zig-zag of quasi-isomorphisms 
$\varphi\mc \Ea \xleftarrow{\eta} \cS \to \Et$. 
Note that there is an obvious bijection between the object sets $\Ba$
and $\Bt$ which is the first condition required for such a
quasi-isomorphism of rings with many objects.  In this section we label 
the elements of both these sets by the corresponding subgroup of $\bT$.     
Recall that $\Ea$ is generated as an algebra by the elements $f_0^H$ and
$g_{2n +1}^H$ with the relations forced by $\Ea(H,K)$ being trivial,
$\Ea(H,{\bT})$ being  
trivial above degree zero and $\Ea(H,H)$ being trivial above degree one.    
Since $\Et(H,H)$ does not necessarily share this property, we must use
the auxiliary ring 
$\cS$ which is basically constructed from $\Ea$ by removing the relations 
forced by $\Ea(H,H)$ being trivial above degree one. 

\begin{proposition}\label{prop-6.2}
There is a differential graded ring $\cS$ with objects $\Ba$ and 
quasi-isomorphisms 
$ \Ea \xleftarrow{\eta} \cS \xrightarrow{\varphi} \Et$ of differential graded rings which induce the
obvious bijections $\Ba \xleftarrow{\id} \Ba \to \Bt$.
\end{proposition}

Combining Propositions~\ref{prop-simple-alg},~\ref{prop-simple-top} and 
\ref{prop-6.2} and
Corollaries~\ref{cor-Gabriel} and~\ref{cor-Gabriel-t} gives our main
result.

\begin{corollary}\label{cor-main}
There is a zig-zag of Quillen equivalences between rational $\bT$-equivariant
spectra and differential graded objects in $\cA$.
\[ \bT\mbox{-spectra} \simeq_Q \dga \]
\end{corollary}

We first define the intermediary ring $\cS$ and the map $\eta\mc \cS \to \Ea$. 
As mentioned above, $\cS(H,H)$ is not required to be trivial above degree one;
so we remove the relation $g_{2k+1}^H f_0^H=0$.  
Then we must add generators $\bb_{i,j}^H$ to ensure that these products are 
still trivial in homology. 
 
\begin{definition}\label{def-S}
{\em Define the differential graded ring $\cS$ with objects $\Ba$ as
the ring generated by elements 
\begin{enumerate}
\item $f_0^H \in \cS(H,\bT)_0$ with $d(f_0^H)=0$, 
\item $g_{2k+1}^H \in 
\cS(\bT,H)_{2k+1}$  for $k \geq 0$ with $d(g_1^H) = 0$ and   
$d(g_{2k+1}^H) = g_1^H f_0^H g_{2k-1}^H$ for $k > 0$  and
\item $\bb_{2m+1, 2n}^H \in \cS(H,H)_{2n}$
for $m \geq 1$ and $n \geq m +1$ 
with $d(\bb_{2m+1, 2m+2}^H) = g_{2m+1}^Hf_0^H$ and $d(\bb_{2m+1,2n}^H) =
\bb_{2m+1, 2n-2}^H g_1^H f_0^H$ for $n > m +1$  
\end {enumerate}
with relations generated by 
$ x^H y^K =0$ for any elements with $H \neq K$, 
$f_0^H \bb_{2m+1, 2n}^H =0$ and $f_0^H g_{2k+1}^H f_0^H =0$.  
The relations are equivalent to setting $\cS(H, \bT) \iso \mQ[0]$
and $\cS(H,K)=0$.  
 
As mentioned in Definition~\ref{def-Ea}, for ease of notation in this
section we denote $\tg_1$ in $\Ea(\bT, H)$ by $g_1$.  The element
$g_1$ in $\EBCa$ does not appear anywhere in this section.
Define $\eta\mc \cS \to \Ea$ as the ring homomorphism determined by $\eta(f_0^H) = f_0^H$,
$\eta(g_{2k+1}^H)= g_{2k+1}^H$ and $\eta(\bb_{2m+1, 2n}^H)=0$. 
}\end{definition}

\begin{proposition}\label{6.4}
The map $\eta \mc \cS \to \Ea$ is a quasi-isomorphism of differential
graded rings.
\end{proposition}

\begin{proof}
For ease of notation we delete the superscripts $H$ and implicitly apply
these steps for each finite subgroup H.  Because there is no interaction
between different finite subgroups, this should not cause any confusion.  
Since $\eta$ takes the relations in $\cS$ to zero in $\Ea$, it is
a ring homomorphism.  It also commutes with the differential since
in $\Ea$ we have that $d(g_{2k+1}) = g_{2k} = g_1 f_0 g_{2k-1}$ and 
$g_{2k+1}f_0 = 0$ for $k > 0$.  
The maps $\cS(H,K) \xrightarrow{\eta} \Ea(H,K)$, 
$\cS(H,\bT) \xrightarrow{\eta} \Ea(H,\bT)$
and $\cS(\bT,\bT) \xrightarrow{\eta} \Ea(\bT,\bT)$ are isomorphisms.

We next show that $H_*\cS(H,H)\iso\mQ[0] \osum \mQ[1]$, generated by 
$[\id_H]$ and $[g_1 f_0]$.  Below is a table of generating elements in degree
zero through seven in $\cS(H,H)$ and $\cS(\bT,H)$.     
\[
\begin{tabular}{c|c|c|c|c|c|c|c|c}
&0&1&2&3&4&5&6&7 \\ \hline
$\cS(H,H)$& $id_H$&$g_1f_0$& &$ g_3f_0$&$\bb_{3,4}$ &$\bb_{3,4}g_1f_0$&$\bb_{3,6}$&$\bb_{3,4}g_3f_0,\  g_7f_0$ \\ 
&& &&&&$g_5f_0$& $\bb_{5,6}$&$ \bb_{5,6}g_1f_0,\  \bb_{3,6}g_1f_0$ \\ \hline
$\cS(\bT, H) $& &$g_1$& $g_1f_0g_1$&$ g_3$&$ g_1f_0g_3$&$ g_5$&$g_1f_0g_5,\  g_3f_0g_3$& $g_7,\  \bb_{3,4}g_3$ \\ 
& && &&$g_3f_0g_1$ &$\bb_{3,4}g_1$&$g_5f_0g_1, \ \bb_{3,4}g_1f_0g_1$&$\bb_{5,6}g_1,\  \bb_{3,6}g_1$ 
\end{tabular}
\]
These elements have been arranged so that the differential  
is either zero or takes an element to the element in the corresponding spot
one degree below.  Here, $d_n\mc \cS_{n} \to \cS_{n-1}$ is non-zero
for $n=4, 6$ on $\cS(H,H)$ and for $n=3, 5, 7$ on $\cS(\bT,H)$.
In general, because of the relations in $\cS$ any element in $\cS(H,H)$ can
be written as a sum of elements of the form $w$ or $w g_{2k+1} f_0$ 
where $w$ is a word in 
the elements of the set $\{ \bb_{2m+1, 2n} \}$ (including the empty word
$w=\id_H$.)   Again because of the relations, $dw$ has only one term:
\[
d(w \bb_{2m+1,2m+2})= w g_{2m+1} f_0 \mbox{ and } d(w \bb_{2m+1, 2n}) =
w \bb_{2m+1, 2n-2} g_1 f_0 \mbox{ for } n > m+1.
\]  
Hence if $w$ is not the empty word, then 
$w$ is not a cycle.  The elements $w g_{2k+1}f_0$ are cycles;  
if $k > 0$ then it is also a boundary (of $w \bb_{2k+1,2k+2}$) 
or if $k = 0$ and $w$ is not
$\id_H$ then again it is a boundary since $w = w'\bb_{2m+1, 2n}$ and 
$d(w'\bb_{2m+1, 2n+2}) = w' \bb_{2m+1, 2n} g_1 f_0 = w g_1 f_0$.       
This leaves only $\id_H$ and $g_1 f_0$ as generators of homology. 
Since $m_1=g_1f_0$ in $\Ea$, we see that $\cS(H,H) \xrightarrow{\eta}
\Ea(H,H)$ induces an isomorphism in homology. 

Finally we show that $H_*\cS(\bT,H) \iso \mQ[1]$ generated by $g_1$ by 
showing that $\cS(\bT, H)\iso (\Osum_{I} D^{n_i}) \osum \mQ[1]$ as a rational 
chain complex.   Every element of $\cS(\bT,H)$ can be written as
a sum of elements of the form $w g_i$ or $w g_j f_0 g_k$ where $i,j$ and
$k$ are odd numbers and $w$ is a (possibly empty) word in the $\bb_{2m+1,2n}$
as above.  Except for $g_1$ these module generators fall into certain families 
$\cF_{w, i, j-i}$ indexed by a word, $w$, and two odd integers, $i$ and $j-i$.  
The families
with $i=1$ are exceptional, $\cF_{w, 1, j-1}$ is the set $\{ g_{j+1}, 
g_1 f_0 g_{j-1}\}$ which generates a complex isomorphic to $D^{j+1}$.   
If $i > 1$ and $j-i = 2l +1$, then $\cF_{w, i, j-i}$ generates a complex 
isomorphic to $l + 1$ copies of $D^{|w| + j + 1}$.
This set is defined by
\[ \cF_{w, i, j-i} 
= \{w \bb_{i,k} g_{j -k +1} \mbox{ for } k \mbox{ even and } i + 1 \leq k \leq j, 
\]\[
w \bb_{i,k} g_1 f_0 g_{j-k-1}  
\mbox{ for } k \mbox{ even and }i+ 1 \leq k \leq j -2,\mbox{ and } w g_i f_0 g_{j-i} \}. 
\]
Since every module generator except $g_1$ appears in one and only one of these families,
this gives the splitting as claimed.  Since $\eta(g_1)=g_1$, it follows
that $\eta\mc \cS(\bT,H) \to \Ea(\bT, H)$ is a quasi-isomorphism.  
\end{proof}

To finish the proof of Proposition~\ref{prop-6.2} we must construct
$\varphi \mc \cS \to \Et$.  This construction uses Toda brackets 
in $\Ho(\Mod\Et)$ which we calculate using 
Massey products.
For $\cR$ a differential graded ring with objects $\cG= \{G_i\}_{i\in I}$,
a cycle $a \in \cR(G_j, G_i)_{p}$ induces a map of right
$\cR$-modules, $a \mc \Sigma^p F_{G_j}^{\cR} \to F^{\cR}_{G_i}$, 
by left multiplication (or 
post-composition) where $F^{\cR}_{G_i}(-)=\cR(-, G_i)$ is the free right 
$\cR$-module represented by $G_i$.     
The associated homotopy class of maps $[a]$ is determined by
the homology class represented by $a$ in the homology of $\cR$. 
This translates between Massey products, which are
homology classes in $\cR$, and Toda brackets, which are homotopy classes
of right $\cR$-module maps.  Notice, if $a \in \cR(G_j, G_i)_p$ and
$b \in \cR(G_i, G_k)_q$ are cycles then extra suspensions are required to 
compose these elements; the associated composition is $b(\Sigma^{q}a) \mc
\Sigma^{p+q}G_j \to \Sigma^q G_i \to G_k$.

Starting with the Massey products for
$\Ea$ we construct elements of the Toda brackets for $\Ho(\Mod\Ea)$
in Proposition~\ref{prop-TB}.
See Definition~\ref{def-Toda} for the definition of higher Toda brackets.
We denote the Toda brackets by $\langle-\rangle_T$, to differentiate them from
the Massey products.  
Since Toda brackets are invariant under exact equivalences by 
Theorem~\ref{thm-tb-eq}, these Toda brackets agree with those for
$\Ho(\Mod\Et)$.  This follows because   
Quillen equivalences induce triangulated equivalences~\cite[I.4.3]{Q},
$\Mod\Et$ is Quillen equivalent to $\bT$-spectra by
Proposition~\ref{prop-simple-top} and Corollary~\ref{cor-Gabriel-t},
$\Mod\Ea$ is Quillen equivalent to $\dga$ by
Proposition~\ref{prop-simple-alg} and Corollary~\ref{cor-Gabriel}
and $\Ho(\bT\mbox{-spectra})$ is triangulated equivalent to $\cD(\cA)$
by~\cite[5.6.1]{Greenlees}.
The Toda brackets in $\Ho(\Mod\Et)$ then give us the structure we
need to construct $\varphi$.

\begin{proposition}\label{prop-TB}
The following Toda brackets in $\Ho(\Mod\Ea)$ agree up to sign
with the Massey products in $H_*\Ea$ from 
Proposition~\ref{prop-homology-brackets}.  
\[
\langle [f_0^H], [m_1^H], \Sigma [m_1^H], \dots, \Sigma^{n-1} [m_1^H], \Sigma^n [g_1^H] 
\rangle_T = \{   [f_0^H g_{2n+1}^H] \} = \{[i_{2n+1}^H]\}
\]
\end{proposition}

The sign difference here corresponds to the different sign conventions 
used for Massey products and Toda brackets.

\begin{proof}
For ease of notation we fix a finite subgroup $H$ and 
delete the superscripts $H$. 
We first consider the three-fold Toda bracket. 
Starting with the defining system for the Massey product 
$\langle [f_0], [m_1],[g_1] 
\rangle$ we construct an element of the Toda bracket 
$\langle [f_0], [m_1], \Sigma [g_1] \rangle_T$.  
By Definition~\ref{def-Toda} such a Toda bracket
is the set of homotopy classes of  compositions 
$\Sigma^3 F_{\bT}^{\Ea} \varr{\beta} X 
\varr{\alpha}
F_{\bT}^{\Ea}$ where the homotopy type of $X$ is a 2-filtered object in 
$\{[m_1]\}$; see Definition~\ref{def-filtered}.  
The cofiber of $m_1$ is such an object. 
Explicitly, set $X= \Sigma^2 F_{H}^{\Ea} \ssm_{m_1} F_{H}^{\Ea}$ with 
underlying
graded object $\Sigma^2 F_{H}^{\Ea} \osum F_H^{\Ea}$ and with differential
given by $d(x,y) = (dx, dy - m_1 x)$.   Then the defining system found
in Proposition~\ref{prop-homology-brackets} 
determines the other two maps as well, $\alpha(x,y) = f_0 y$
and $\beta(z)= (g_1 z, g_3 z)$.  Here we are using the sign convention that
$d_{\Sigma C} = -d_C$.  One can check that $[\sigma_X \beta] = [g_1] $ and
$[\alpha \sigma_X']=[f_0]$.  
Due to the sparseness of elements in $\Ea$, one can check that $\alpha, 
\beta$ are the unique
maps satisfying the definition of the Toda bracket. Since $\Sigma^3 
F_{\bT}^{\Ea}, X$ and $F^{\Ea}_{\bT}$ are cofibrant and fibrant 
$\Ea$-modules
these are also the unique such homotopy classes. 
Hence the Toda bracket 
is the one element set $\{ [f_0 g_3] = [i_3] \}$.  

As with the three-fold Massey product, the higher Massey products
also agree up to sign with the higher Toda brackets.  
First, define the iterated cofiber
\[X_n =
\Sigma^{2n} F_H^{\Ea} \ssm_{m_1} \Sigma^{2n-2} F_H^{\Ea} \ssm_{m_1} \dots
\ssm_{m_1} F_H^{\Ea}
\]
with differential $d(x_1, x_2, \dots, x_{n+1})=
(dx_1, dx_2 - m_1 x_1, \dots, dx_{n+1}- m_1 x_n)$.  Note here $d^2=0$ 
because $m_1^2 =0$.  By Lemma~\ref{lem-cohen},
$X_n$ is an element of $\{ [m_1], \Sigma [m_1], \dots, \Sigma^{n-1} [m_1]\}$.
Then the associated
defining system determines the maps $\Sigma^{2n+1} F_{\bT}^{\Ea}
\varr{\beta} X_n \varr{\alpha} F_{\bT}^{\Ea}$ as $\beta(z)=(g_1 z, g_3 z, \dots, g_{2n+1}z)$
and $\alpha(x_1, \dots, x_{n+1})=f_0 x_{n+1}$.  This shows that the composite,
$[f_0 g_{2n+1}]=[i_{2n+1}]$ is in the Toda bracket $\langle [f_0], [m_1], 
\Sigma [m_1], \dots, \Sigma^{n-1} [m_1], \Sigma^n [g_1] \rangle_T$.  One can
check that any other defining system would produce a Toda bracket and
vice versa, but again the sparseness of elements in $\Ea$ forces this to 
be 
the only element in the Toda bracket as well as the only defining system.  
\end{proof}
\begin{proof}[Proof of Proposition~\ref{prop-6.2}]
We are left with constructing $\varphi\mc \cS \to \Et$.  We proceed by
induction on the degree.    
We exploit the fact that Toda brackets in $\Ho(\Mod\Et)$ and
$\Ho(\Mod\Ea)$ agree to construct $\varphi$.  At each stage certain
products must be non-trivial because of the non-trivial Toda brackets.
This then provides non-trivial targets for $\varphi$ of the elements 
$g_{2k+1}^H$ in $\cS$.  
The products $g_{2k+1}^H f_0^H$ are trivial in homology (for $k > 0$)
so they must be boundaries.  This produces targets for $\varphi$ of the
elements $\bb_{2k+1, 2l}^H$ in $\cS$.   We can then extend 
$\varphi$ to a ring homomorphism because the images of the relations in 
$\cS$ between these generators hold in $\Et$ as well since  $\Et(H, \bT)$ is 
concentrated in degree zero.   
By construction $\varphi$ is then a map of differential graded objects and
we show inductively that it is also a quasi-isomorphism.

Again, we mostly delete the superscripts $H$ and 
implicitly apply these steps for each finite subgroup. 
Since $H_*\Et$ is isomorphic to $H_*\E(\Bt)$, 
we use the notation for elements in $H_*\Et$ from Proposition~\ref{prop-homology-top}.  To further ease notation, if $y$ is the boundary of $x$, 
i.e. $dx =y$, then we say that $x$ is a {\em primitive} of $y$. 

The images of $\id_H$ and $\id_{\bT}$ in $\cS$ are forced to be
the respective identity elements in $\Et$. 
We next determine the images of $f_0$ and $g_1$. 
Let $f_0'$ be an element in $\Et(H,\bT)_0$ such that $[f_0'] = \wh{f_0}$
in $H_0\Et(H,\bT)$. 
Let $g_1'$ be a cycle such that $[g_1']=\wh{g_1}$ in $H_1\Et(\bT, H)$.  
Then set $\varphi(f_0)=f_0'$ and $\varphi(g_1)=g_1'$.  
Since $\Et(H, \bT)$ is concentrated in degree zero by 
Proposition~\ref{prop-simple-top}, $f_0' g_1' f_0'$ is trivial.
Thus the relations between $f_0$ and $g_1$ in $\cS$ also hold in $\Et$
so we may extend $\varphi$
to products of $g_1$ and $f_0$ using the product structure in $\Et$.
Since $[g_1'][f_0']= \wh{m_1}$ and $[f_0'] [g_1']=\wh{i_1}$ are non-trivial in 
$H_*\Et$ by Propostions~\ref{prop-homology-top} and~\ref{prop-simple-top}, the 
images $\varphi(g_1 f_0)$ and $\varphi(f_0 g_1)$ must be non-trivial cycles. 
Hence $\varphi$ induces a quasi-isomorphism in degrees zero and one.    

Before turning to the general induction step, we consider degrees two and
three.  It may be useful to refer to the tabulation of elements of
$\cS(H,H)$ and $\cS(\bT, H)$ given in the proof of Proposition~\ref{6.4}.
There we also showed the other components of $\cS$ are isomorphic to $\Ea$,
see Definition~\ref{def-Ea}; as vector spaces $\cS(H,K)= 0$, $\cS(H, \bT)$
is generated by $f_0$ and $\cS(\bT, \bT)$ is generated by the elements
$\id_{\bT}, f_0 g_1, f_0 g_3, \dots$. 
As mentioned above, Toda brackets agree for $\Ho(\Mod\Ea)$ and $\Ho(\Mod\Et)$
by Theorem~\ref{thm-tb-eq}.  Since $\langle[f_0], [m_1], \Sigma[g_1]\rangle_T
=\{[i_{2n+1}]\}$ by Proposition~\ref{prop-TB} in $\Ho(\Mod\Ea)$, 
there must be maps $\Sigma^3 F_{\bT}^{\Et} \varr{\beta'} X' \varr{\alpha'} 
F_{\bT}^{\Et}$ such that $[\alpha'\beta'] = \wh{i_3}$ is in  
$\langle \wh{f_0}, \wh{m_1}, \Sigma \wh{g_1} \rangle_T$ in $\Ho(\Mod\Et)$.  
Here $X'$ must be a 2-filtered object in $\{\wh{m_1}\}$, so define
$X'$ as the cofiber of $g_1' f_0'$,
similar to $X$ in the proof of Proposition~\ref{prop-TB} above.  
Then since $\Et(H, \bT)_n$ is trivial for $n > 0$
the map $\alpha'$ must be $\alpha'(x,y) = f_0' y$.
If $g_1' f_0' g_1'=0$ then one could set $\beta'(z) = (g_1' z, 0)$.  Then the
composite $[\alpha'\beta']$ would be trivial.  Since this Toda bracket does not
contain the trivial map, this means $g_1' f_0' g_1'$ must be non-trivial. 
Since $H_2\Et(\bT,H)=0$ it must also be a boundary, 
let $g_3'$ be one of its primitives in $\Et(\bT,H)_3$. 
Then set $\beta'(z) = (g_1' z, g_3' z)$.
The maps $\alpha'$ and $\beta'$ satisfy the definition of the Toda bracket, so 
by Proposition~\ref{prop-TB} we must have $[\alpha'\beta']=[f_0' g_3']=\wh{i_3}$.   
Set $\varphi(g_3)=g_3'$.    

We next inductively determine the images of $\bb_{3,2n}$ with $n \geq 2$.  Since 
$f_0' g_1' f_0'$ is in the trivial group $\Et(H, \bT)_1$, $d(g_3'f_0')=
(dg_3')f_0'=g_1'f_0'g_1'f_0'=0$ and     
the element $g_3' f_0'$ is a cycle.  It is also a boundary since 
$H_3\Et(H,H) =0$ by Proposition~\ref{prop-homology-top} and~\ref{prop-simple-top}; let $\bb_{3,4}'$ be one of its primitives 
and set $\varphi(\bb_{3,4})=\bb_{3,4}'$.  Assume that 
elements $\bb_{3,2k}'$ have been chosen which are
primitives of 
$\bb_{3,2k-2}'g_1'f_0'$ for $2 < k < n$.  
Again, since $f_0' g_1' f_0'$ is trivial and $d(\bb_{3, 2n-2}'g_1'f_0')=
\bb_{3, 2n-4}'g_1'f_0'g_1'f_0'$  the element $\bb_{3,2n-2}'g_1'f_0'$ is a cycle.
Since $H_{2n-1}\Et(H,H)=0$ for $n>1$ by Proposition~\ref{prop-homology-top} and~\ref{prop-simple-top}, 
it is also a boundary.  Let $\bb_{3,2n}'$
be one of its primitives
and set $\varphi(\bb_{3,2n})=\bb_{3,2n}'$.  

We can now extend $\varphi$ to the subring of $\cS$ generated by $f_0, g_1,
g_3$ and $\bb_{3,2n}$ for $n \geq 2$.  This is possible because the relations
among these elements in $\cS$ are also satisfied in $\Et$.  Specifically, since
$\Et(H, \bT)$ is concentrated in degree zero $f_0'\bb_{3,2n}'=0,
f_0'g_1'f_0'=0$ and $f_0'g_3'f_0'=0$.  By construction $\varphi$ also
commutes with the differential on this subring.  The homology
in degrees two and three of $\Et$ is generated by    
the elements $[{f_0^H}' {g_3^H}'] = \wh{i_3^H}$ in 
$H_3\Et= H_3\Et(\bT,\bT)$.   In Proposition~\ref{6.4} we showed $\eta\mc\cS
\to \Ea$ is a quasi-isomorphism, so Propostion~\ref{def-Ea} shows that
the elements $f_0^Hg_3^H$ generate the homology in degrees two and
three for $\cS$ as well.
Thus, $\varphi$ induces a quasi-isomorphism through degree 3.

Assume by induction that $\varphi(g_{2k+1}) = g_{2k+1}'$ and 
$\varphi(\bb_{2k+1, 2l})= \bb_{2k+1, 2l}'$ for $k \leq n$ and $l \geq k+1 $,  
and that the extension of $\varphi$ to the subring generated by these elements  
induces a quasi-isomorphism through degree $2n+1$.  
Since $g_1' f_0'g_1'f_0' = 0$, an $\Et$-module 
\[
X_n' \in \{\wh{m_1}, \Sigma \wh{m_1}, \dots, \Sigma^{n-1}\wh{m_1}\}
\] 
exists, and can be built similarly to $X_n$ in Proposition~\ref{prop-TB}.  
Since $\Et(H, \bT)_n=0$
for $n > 0$, the map $\alpha' \mc X_n' \to F_{\bT}^{\Et}$ is determined as
$\alpha'(x_1, \dots, x_{n+1})= f_0' x_{n+1}$.  As with the three-fold Toda 
bracket
if $g_1' f_0' g_{2n+1}'$ were trivial then we could set $\beta'(z)=
(g_1' z, g_3' z, \dots, g_{2n+1}' z, 0)$.  But then the composite $\alpha'\beta'$ 
would be trivial, which is not allowed because the $n+2$-fold
Toda bracket from Proposition~\ref{prop-TB} does not contain the
trivial map.  Hence $g_1' f_0' g_{2n+1}'$ is
non-trivial.  Since $H_{2n+2}\Et(\bT, H)=0$, it must also be a boundary; let  
$g_{2n+3}' \in \Et(\bT, H)_{2n+3}$ be one of its primitives
and set $\beta'(z) = (g_1' z, g_3' z, \dots, g_{2n+3}'z)
$.  Then $\alpha'$ and $\beta'$ satisfy the definition of the Toda Bracket so 
\[
[\alpha'\beta']=[f_0' g_{2n+3}'] \in \langle \wh{f_0}, \wh{m_1}, 
\Sigma \wh{m_1}, \dots, \Sigma^{n-1} \wh{m_1}, \Sigma^n \wh{g_1} \rangle_T
= \{ \wh{i_{2n+3}}\}.
\]  
Set $\varphi(g_{2n+3})= g_{2n+3}'$. 

Arguments similar to those for $\bb_{3,2k}$ for $k \geq 2$ also apply to the 
elements $\bb_{2n+3,2k}$ for $k\geq n+2$, $n > 0$.  
The element $g_{2n+3}' f_0'$ is a cycle and a boundary since 
$H_{2n+3}\Et(H,H) =0$; let $\bb_{2n+3,2n+4}'$ be one of its primitives 
and set $\varphi(\bb_{2n+3,2n+4})=\bb_{2n+3,2n+4}'$.  Assume that
elements $\bb_{2n+3,2k}'$ have been chosen which are primitives of
$\bb_{2n+3,2k-2}' g_1'f_0'$ for $n+2<k < l$.  Since $f_0'g_1'f_0'$ is trivial and $H_{2n-1}\Et(H,H)=0$, 
the element $\bb_{2n+3,2l-2}'g_1'f_0'$ is a cycle and a boundary.  Let 
$\bb_{2n+3,2l}'$
be one of its primitives 
and set $\varphi(\bb_{2n+3,2l})=\bb_{2n+3,2l}'$. 

We can now extend $\varphi$ to the subring of $\cS$ generated by $f_0$, $g_{2k+1}$
and $\bb_{2k+1,2l}$ for $k \leq n+1$ and $l \geq k+1$; again, this is possible 
because $\Et(H, \bT)$ is concentrated in degree zero, so the relations
hold.  Also, by construction $\varphi$ commutes with the differential on
this subring. 
Since the elements $[{f_0^H}' {g_{2n+3}^H}']=\wh{i_{2n+3}^H}$ generate 
$H_{2n+3}\Et=H_{2n+3}\Et(\bT,\bT)$ and the elements $f_0^H g_{2n+3}^H$
generate $H_{2n+3}\cS\iso H_{2n+3}\Ea \iso H_{2n+3}\Ea(\bT,\bT)$,
we see that $\varphi$ induces a quasi-isomorphism through degree $2n+3$.
Hence by induction we have constructed the quasi-isomorphism $\varphi
\mc \cS \to \Et$. 
\end{proof}

\begin{appendix}
\section{Higher Toda brackets and Massey products}\label{sec-Toda}

In this section we show that an exact equivalence of 
triangulated categories preserves higher Toda brackets.  
We first make slight modifications to the definitions of higher 
Toda brackets 
from~\cite[2]{cohen} so that they apply to an arbitrary triangulated 
category $\cT$.  We have kept the same order for the maps in the
Toda bracket as in~\cite{cohen}, but we have reversed the numbering
to match the Massey products numbering.

\begin{definition}\label{def-filtered}
{\em An {\em $n$-filtered object} $X$ is an object $X$ in $\cT$ 
with maps $*\iso F_0 X \varr{i_0} F_1 X \varr{i_1} \cdots \to F_{n-1} X 
\varr{i_{n-1}} F_nX \iso X$.  Let $F_{k+1}X/F_kX$ denote the cofiber of
$i_k$ so that there are triangles, 
$F_{k} X \varr{i_k} F_{k+1} X \varr{\pi_{k+1}} F_{k+1}X/F_kX \varr{\delta_k}
\Sigma F_{k} X$.
Given a composable sequence of maps, 
$A_{n-1} \varr{f_{n-1}} A_{n-2} \varr{f_{n-2}} \cdots A_1\varr{f_{1}} 
A_0$, an $n$-filtered object $X$ {\em is an element of 
$\{ f_{1}, \dots, f_{n-1}\}$} if and only if $\Sigma^k A_{k} \iso F_{k+1}X/F_kX$
and via these isomorphisms $\Sigma^k{f_{k}}$ is isomorphic to the composite
$F_{k+1}X/F_kX \varr{\delta_k} \Sigma F_{k} X \varr{\Sigma \pi_k} \Sigma
(F_kX /F_{k-1}X)$.
For $X \in \{f_{1}, \dots, f_{n-1}\}$, define $\sigma_X'$ as the composite
$A_0 \iso F_1X \to X$ and define $\sigma_X$ as the composite $X\iso F_nX 
\varr{\pi_n} F_n X/F_{n-1} X \iso \Sigma^{n-1} A_{n-1}$.
}
\end{definition}

\begin{definition}\label{def-Toda}
{\em Define the $n$-fold Toda bracket $\langle f_1, f_2, \dots, 
f_n\rangle_T$ for a sequence of maps
\[ A_n \varr{f_n}A_{n-1} \varr{f_{n-1}} \cdots \varr{f_{2}} A_1 
\varr{f_1} A_{0} \]
as the set of all maps $\theta \in \cT(\Sigma^{n-2} A_n, A_{0})$ such that
there is some $(n-1)$-filtered object $X \in \{f_{2},\dots, f_{n-1}\}$ and maps 
$\beta\mc \Sigma^{n-2}A_n\to X$ and $\alpha \mc X \to A_{0}$ with $\Sigma^{n-2} f_n
\iso \sigma_X \beta$, $f_1 \iso \alpha\sigma_X'$ and $\theta \iso \alpha\beta$. 
\[
\xymatrix{
{}&{A_1} \ar[d]_{\sigma_X'}\ar[dr]^{f_1}& {}\\
{\Sigma^{n-2}A_n} \ar[r]^{\beta} \ar[dr]_{\Sigma^{n-2}f_n} & X \ar[d]^{\sigma_X}
\ar[r]^{\alpha} & {A_0}\\
&{\Sigma^{n-2}A_{n-1}}
}
\]
}\end{definition}

With these definitions it is easy to see that Toda brackets are preserved by
exact equivalences of triangulated categories.

\begin{theorem}\label{thm-tb-eq}
Given an exact equivalence of triangulated categories $\psi \mc \cT \to \cT'$,
then $\theta \in \langle f_1, \dots, f_n \rangle_T$ if and only if 
$\psi(\theta) \in \langle \psi(f_1), \dots, \psi(f_n) \rangle_T$.
\end{theorem}

\begin{proof}
Assume given $\theta \in \langle f_1, \dots, f_n \rangle_T$.  Then $\theta = \alpha \beta$
with $\beta\mc \Sigma^{n-2}A_n \to X$ and $\alpha \mc X \to A_{0}$ with
$X \in \{f_{2}, \dots, f_{n-1}\}$.  Since $\psi$ is exact, 
$\psi(X)$, with filtrations $\psi(F_kX)$, is an element of $\{\psi(f_{2},
\dots, \psi(f_{n-1})\}$.   Hence $\psi(\theta) = \psi(\alpha) \psi(\beta)$ is an
element of $\langle \psi(f_{1}), \dots, \psi(f_n)\rangle_T$.  For the other
direction, apply the same arguments to the inverse equivalence $\psi^{-1}$.
\end{proof}

The following manipulations are useful for constructing filtered objects. 
Let $C_{\alpha}$ denote the cofiber of the map $\alpha$.

\begin{lemma}\label{lem-cohen}
Assume given a composable sequence of maps $f_i\mc A_i \to A_{i-1}$.
\begin{enumerate}
\item If $X \in \{f_{2}, \dots, f_k\}$ and $\alpha \mc X \to A_{0}$,
then $C_{\alpha} \in \{ \alpha \sigma_X', f_{2}, \dots, f_k\}$.
\item If $Y \in \{ f_1, \dots, f_{k-1}\}$ and $\alpha\mc \Sigma^{k-1} A_k \to Y$,
then $C_{\alpha} \in \{f_1, \dots, f_{k-1},  \Sigma^{-k+1} \sigma_Y \alpha \}$.
\item If $X \in \{f_1, \dots, f_k\}$, then $F_kX \in \{f_1, \dots, f_{k-1}\}$
and $\Sigma^{-1}(X/F_1X) \in \{f_2, \dots, f_{k}\}$.
\end{enumerate}
\end{lemma}

\begin{proof} The first two statements follow 
by setting $n=1$ and $m=1$ respectively in~\cite[Proposition 2.3]{cohen}.
The third statement follows from~\cite[Proposition 2.2]{cohen} and its analogue.\end{proof}

Often definitions of Toda brackets require various vanishing conditions.
The following proposition shows that these are implicit in the above definition.

\begin{proposition}\label{prop-zero}
There is an object $X$ in $\{f_{1}, \dots, f_{n-1}\}$ if and only if
$0 \in \langle f_{1}, \dots, f_{n-1}\rangle_T$.
\end{proposition}

\begin{proof}
First assume $X \in \{f_{1}, \dots, f_{n-1}\}$.   By Lemma~\ref{lem-cohen}, 
$\Sigma^{-1}(F_{n-1}X /F_1X) \in \{f_{2}, \dots, f_{n-2} \}$.  Let $\beta$ 
be the composite $\Sigma^{n-3}A_{n-1} \iso \Sigma^{-2}(F_nX/F_{n-1}X) \to 
\Sigma^{-1}F_{n-1} X \to \Sigma^{-1}(F_{n-1}X/ F_1 X)$.  Let $\alpha$ be 
the 
map $\Sigma^{-1}(F_{n-1}X/ F_1 X) \to F_1 X \iso A_0$.  Since $\alpha$ is 
the
cofiber of the map $\Sigma^{-1}F_{n-1} X \to \Sigma^{-1}(F_{n-1}X/ F_1 X)$,
the composite $\alpha \beta$ is trivial.  So $\alpha \beta = 0 \in 
\langle f_{1}, \dots, f_{n-1} \rangle_T$.    

For the other direction, assume $\theta= 0 \in \langle f_{1}, \dots,
f_{n-1} \rangle_T$.  Then $\theta = \alpha \beta$ with $\beta \mc \Sigma^{n-3} A_{n-1} 
\to Z$
and $\alpha\mc Z \to A_0$ with $Z \in \{f_{2}, \dots, f_{n-2} \}$.  Since $\alpha \beta$
is trivial, taking cofibers gives a map 
$\gamma \mc \Sigma^{n-2} A_{n-1} \to 
C_{\alpha}$.  By Lemma~\ref{lem-cohen}, $C_\alpha \in \{f_{1}, \dots, f_{n-2}\}$ and
$C_{\gamma} \in \{f_{1}, \dots, f_{n-1}\}$.
\end{proof}

\end{appendix}

\end{document}